\documentclass{amsart} 
\usepackage{amsmath,amssymb,amsxtra}
\usepackage[dvips]{graphics}
\usepackage{enumerate}
\usepackage{array}
\usepackage[matrix,arrow]{xy} 

\newcommand{\finiteG}{{\textsf{G}}}
\newcommand{\finiteH}{{\textsf{H}}}
\newcommand{\finiteL}{{\textsf{L}}}
\newcommand{\finiteC}{{\textsf{C}}}
\newcommand{\finiteS}{{\textsf{S}}}
\newcommand{\finiteT}{{\textsf{T}}}

\newcommand{\finiteM}{{\textsf{M}}}
\newcommand{\finiteP}{{\textsf{P}}}
\newcommand{\finiteB}{{\textsf{B}}}
\newcommand{\ind}{\textup{ind}}
\newcommand{\Ind}{\textup{Ind}}
\newcommand{\infl}{\textup{infl}}

\newcommand{\C}{{\mathbb{C}}}

\newcommand{\Hom}{\textup{Hom}}
\newcommand{\OK}{\mathcal{O}}
\newcommand{\N}{{\mathcal{N}}}
\newcommand{\St}{\textup{St}}
\newcommand{\A}{\mathcal{A}}
\newcommand{\tA}{\widetilde{\A}}
\newcommand{\B}{{\mathcal{B}}}
\newcommand{\D}{\mathcal{D}}
\newcommand{\F}{{\mathcal{F}}}

\newcommand{\reg}{^{\textup{reg}}}
\newcommand{\ereg}{^{\varepsilon\textup{-reg}}}
\let\det\relax
\newcommand{\det}{\textup{det}}

\newcommand{\epsgen}{\Gamma} 
\newcommand{\Gal}{\textup{Gal}}
\newcommand{\phm}{\phantom{-}}
\newcommand{\phinv}{\phantom{{}^{-1}}}
\newcommand{\ptp}{{\widetilde\pi'}}

\newcommand{\reducible}[1]{(\ref{sec:depth-zero-G}PS--#1)}
\newcommand{\depthzerosc}[1]{(\ref{sec:depth-zero-G}SC--#1)}
\newcommand{\cartan}[1]{(\ref{sec:cartan-G}--#1)}

\newcommand{\GL}{\text{GL}}
\newcommand{\PGL}{\text{PGL}}
\newcommand{\UU}{U}

\newtheorem{thm}{Theorem}[section] 

\newtheorem{prop}[thm]{Proposition} 
\newtheorem{lem}[thm]{Lemma} 
\newtheorem{cor}[thm]{Corollary} 
\theoremstyle{definition}

\theoremstyle{remark} 
\newtheorem{rem}[thm]{Remark} 

\numberwithin{equation}{subsection}

\newlength\labeldistance
\begin{document}

\title
{Depth-zero base change for unramified $U(2,1)$}
\author{Jeffrey D. Adler}
\thanks{The first-named author was partially supported
by the National Security Agency (\#MDA904-02-1-0020).}
\address{Department of Theoretical and Applied Mathematics \\
The University of Akron \\
Akron, OH  44325-4002}
\email{adler@uakron.edu}
\author{Joshua M. Lansky} 
\address{Department of Mathematics and Statistics\\
American University\\
Washington, DC 20016}
\email{lansky@american.edu} 


\begin{abstract}
We give an explicit description of $L$-packets and quadratic base
change for depth-zero representations of unramified unitary groups
in two and three variables.  We show that this base change is
compatible with unrefined minimal $K$-types.
\end{abstract}

\maketitle

\section{Introduction}
\label{sec:intro}
Given a finite Galois extension $E/F$ of finite, local, or global fields
and a reductive algebraic $F$-group $\underline{G}$,
``base change'' is, roughly, a (sometimes only conjectural) 
mapping from
representations 
of $G = \underline{G}(F)$
to those of $\underline{G}(E)$.
When $F$ is finite, or when $F$ is local
and $\underline{G}= \GL(n)$,
then this mapping is the Shintani lifting
(as introduced in~\cite{shintani}
and extended in~\cite{kawanaka},~\cite{gyoja},
and~\cite{digne0} for finite groups).

Correspondences like base change that are associated to
the Langlands program can be difficult to describe explicitly,
even in cases where they are known to exist.
Bushnell and Henniart~\cite{%
bushnell-henniart:local-tame-1,%
bushnell-henniart:local-tame-2,%
bushnell-henniart:local-tame-4,%
bushnell-henniart:bc-glp%
}
are remedying this situation for base change for $\GL(n)$ over local fields.
Analogously, Silberger and Zink
\cite{silberger-zink:level-zero-matching}
have made the Abstract Matching
Theorem~\cite{dkv:amt,rogawski:amt,badulescu:thesis}
explicit for depth-zero discrete series representations.

Suppose that $F$ is a $p$-adic field of odd residue characteristic.
If $E/F$ is quadratic, and $\underline{G}$ is a unitary group
in three variables defined with respect to $E/F$, then
Rogawski~\cite{rogawski} has shown that a base change lifting
exists, and has derived some of its properties.
Our goal in this paper is to describe base change explicitly
for depth-zero representations
in the case where $E/F$ is unramified.
Depth-zero base change is particularly interesting
because it should be closely related to base change
for finite groups.  See~\cite{kim-ps:theta-10}
for an exploration of another special case
of this phenomenon.

In order to apply a technical lemma (Cor.~\ref{cor:q-value}), we will
assume that the order $q$ of the residue field $k_F$ of $F$ is at
least $59$.  From the
lemma, character identities can be verified by evaluation at ``very
regular'' elements.  At such elements, character values are
particularly easy to compute.  Without the lemma, the verification of
these identities involves evaluation at more general elements.
Character values at these elements can be computed, but are far more
complicated.

Note that we assume that $F$ has characteristic zero only so that we
can apply results of Rogawski~\cite{rogawski}.  Our calculations apply
equally well if $F$ is a function field of odd residue characteristic.

Conjecturally, one should be able to determine the depth
of the representations in an $L$-packet from the associated
Langlands parameter.  Thus, liftings that arise from the Langlands
correspondence should preserve depth, if depth is normalized correctly.
In particular, depth-zero representations should go to depth-zero
representations.  We assume this throughout for base change
and for endoscopic lifting from $U(1,1)\times U(1)$ to $U(2,1)$.

In \S\ref{sec:prelim}, we present our notation, review the general
notion of Shintani lifting, describe how it applies to the
representations of certain finite reductive subquotients of $G$, and
list all of the representations of $G$ of depth zero.  In
\S\ref{sec:L-packets}, we give an explicit description of the
depth-zero $L$-packets and $A$-packets for $G$.  In
\S\ref{sec:base_change}, we determine the base change lift of each of
these packets.  In \S\ref{sec:K-types} we examine the relationship
between base change and \emph{$K$-types}, as defined by
Bushnell-Kutzko~\cite{bushnell-kutzko:smooth} and as described by
Moy-Prasad~\cite{moy-prasad2} or Morris~\cite{morris:intertwining}.
Recall that a (minimal) $K$-type (or simply a ``type'') of depth zero
is a pair $(G_x,\sigma)$, where $G_x$ is a parahoric subgroup of $G$,
and $\sigma$ is the inflation to $G_x$ of an irreducible cuspidal
representation of the finite reductive quotient $\finiteG_x$ of $G_x$.
Since all of the data in this definition can be lifted in a natural
way to similar data for $\widetilde{G} = \underline{G}(E)$, we have a
natural notion of base change for depth-zero types.  Under the
above assumption on the residue characteristic of the $p$-adic field
$F$, we show that base change for depth-zero types is compatible with
base change for representations (actually, $A$-packets of
representations):
\begin{thm}
\label{thm:main}
Suppose $\Pi$ is a depth-zero $A$-packet for $G$,
let $\tilde\pi$ denote the base-change lift of $\Pi$,
and let $\pi\in\Pi$.
Suppose $(G_x,\infl(\sigma))$ is a type contained in $\pi$.
Then $\tilde\pi$ contains $(\widetilde G_x,\infl(\tilde\sigma))$,
where $\tilde\sigma$ is the base-change lift of $\sigma$
from $\finiteG_x$ to $\widetilde{\finiteG}_x$.
\end{thm}

Note that the pair $(\widetilde{G}_x,\infl(\tilde\sigma))$ contains a
type upon restriction to some parahoric subgroup of $\widetilde{G}_x$.
Thus, either it is itself a type, or it carries more information than
a type.

In \S\ref{sec:induced}, we state a formula for the character
of an induced representation.  The formula itself is not new,
but we need to assert that it holds for representations of groups
that are not necessarily connected.

In order to describe explicit base change for all representations of $U(2,1)$
(not just of depth zero), one needs to understand depth-zero base change
not just for $U(2,1)$ but for unitary groups in two variables as well.
We deal with this briefly in \S\ref{sec:u11}.

We thank Robert Kottwitz, Jonathan Rogawski, A.~Raghuram,
David Pollack, Stephen DeBacker, and a referee for helpful communications.

\begin{table}

\rotatebox{90}{
\begin{tabular}{|l r| l r|}
\multicolumn{2}{c}{$A$-packet$^*$} & \multicolumn{2}{c}{Base change lift} \\
\hline
$\{\Ind_B^G \lambda\}$ ($\Ind_B^G \lambda$ irreducible and
        $\Ind_{\widetilde B}^{\widetilde G} \tilde\lambda$ irreducible)
        & (\S\ref{sub:ps})
        & $\Ind_{\widetilde B}^{\widetilde G} \tilde\lambda$
        & (\S\ref{sub:bc_ps})
\\
\hline
$\{\Ind_B^G \lambda\}$ ($\Ind_B^G \lambda$ irreducible and
        $\Ind_{\widetilde B}^{\widetilde G} \tilde\lambda$ reducible)
        & (\S\ref{sub:ps})
        & $\Ind_{\widetilde P}^{\widetilde G}
        \Bigl(
        ( \lambda_1\tilde\lambda_2|\cdot|_E^{\mp 1/2}\circ \det_{\GL(2)})
                \otimes \tilde\lambda_2
        \Bigr)$
        & (\S\ref{sub:bc_ps})
\\
\hline
$\{\psi\}$ (one-dimensional)
        & (\S\ref{sub:ps})
        & $\tilde\psi$
        & (\S\ref{sub:bc_ps})
\\
\hline
$\{\St_G(\psi)\}$
        & (\S\ref{sub:ps})
        & $\St_{\widetilde G}(\tilde\psi)$
        & (\S\ref{sub:bc_ps})
\\
\hline
$\{\pi_1(\lambda),\pi_2(\lambda)\}$
        & (\S\ref{sub:ps})
        & $\Ind_{\widetilde B}^{\widetilde G} \tilde\lambda$
        & (\S\ref{sub:bc_ps})
\\
\qquad ($\Ind_B^G \lambda = \pi_1(\lambda) \oplus \pi_2(\lambda)$) &&&
\\
\hline
$ \Ind_{G_y}^{G}\sigma$, $\sigma$ cubic cuspidal
        & (\S\ref{sub:stable_sc})
        & $\Ind_{\widetilde Z\widetilde G_y}^{\widetilde G} \tilde\sigma$
        & (\S\ref{sub:bc_stable_sc})
\\
\hline
$\{\pi^2(\lambda),\pi^s(\lambda)\}$
        & (\S\ref{sub:ps}) 
        & $\Ind_{\widetilde P}^{\widetilde G} \St_{\widetilde H}
        \Bigl(
        ( \lambda_1\tilde\lambda_2|\cdot|_E^{\mp 1/2}\circ \det_{\GL(2)})
                \otimes \tilde\lambda_2
        \Bigr)$
        & (Prop.~\ref{prop:bc_card_2})
\\
($\pi^s$ described in Prop.~\ref{prop:card2}) &&&
\\
\hline
$\{\pi^n(\lambda),\pi^s(\lambda)^{(*)}\}$
        & (\S\ref{sub:ps}) 
        & $\Ind_{\widetilde P}^{\widetilde G} 
        \Bigl(
        ( \lambda_1\tilde\lambda_2|\cdot|_E^{\mp 1/2}\circ \det_{\GL(2)})
                \otimes \tilde\lambda_2
        \Bigr)$
        & (Prop.~\ref{prop:bc_card_2})
\\
($\pi^s$ described in Prop.~\ref{prop:card2}) &&&
\\
\hline
$
\begin{array}{@{}lcr}
\bigl\{
&\Ind_{G_y}^G \infl_{\finiteG_y}^{G_y}
(-R_\finiteC^{\finiteG_y} \varphi_1\otimes\varphi_2\otimes\varphi_3), \\
&\Ind_{G_z}^G \infl_{\finiteG_z}^{G_z}
(-R_{\finiteC'}^{\finiteG_z} \varphi_1\otimes\varphi_2\otimes\varphi_3), \\
&\Ind_{G_z}^G \infl_{\finiteG_z}^{G_z}
(-R_{\finiteC'}^{\finiteG_z} \varphi_2\otimes\varphi_3\otimes\varphi_1), \\
&\Ind_{G_z}^G \infl_{\finiteG_z}^{G_z}
(-R_{\finiteC'}^{\finiteG_z} \varphi_3\otimes\varphi_1\otimes\varphi_2) &\bigr\}
\end{array}
$
        & (Prop.~\ref{prop:card4})
        & $\Ind_{\widetilde Z \widetilde G_y}^{\widetilde G}
           \infl_{\widetilde\finiteG_y}^{\widetilde Z \widetilde G_y}
        (-R_{\widetilde\finiteC}^{\widetilde \finiteG_y}
                \tilde\varphi_1\otimes\tilde\varphi_2\otimes\tilde\varphi_3)$
        & (Prop.~\ref{prop:bc_card4})
\\
\qquad($\varphi_i$ distinct) &&&
\\
\hline
\end{tabular}
}
\caption{Depth-zero $A$-packets for $U(2,1)$, and their base change lifts.
To obtain $L$-packets,
omit the representations marked with an asterisk $(*)$.}
\label{table:results}
\end{table}

\section{Preliminaries}
\label{sec:prelim}
\subsection{General notation and facts}
\label{sec:notation-general}
For any nonarchimedean local field $F$, let $\OK_F$ denote
its ring of integers, $\mathfrak{p}_F$ the prime ideal
in $\OK_F$, and $k_F= \OK_F/\mathfrak{p}_F$ the
residue field.
For any abelian extension $E/F$, let
$\omega_{E/F}$ denote the character
of $F^\times$ arising via local class field theory.

We will use underlined letters to denote algebraic groups and will
drop the underlining to indicate the corresponding groups of rational
points.  Given an algebraic $F$-group $\underline{G}$ and a finite
extension $E/F$, let $\widetilde{\underline{G}} =
R_{E/F}(\underline{G})$, where $R_{E/F}$ denotes restriction of
scalars.  Similarly, if $\underline\finiteG$ is a $k_F$-group,
$\widetilde{\underline{\finiteG}}$ will denote
$R_{k_E/k_F}(\underline{\finiteG})$.  Whenever we use this notation, the
extension $E/F$ will either be specifed, or it will be understood from
the context.

For every nonarchimedean local field $F$ and every reductive algebraic
$F$-group $\underline{G}$, one has an associated \emph{extended affine
  building} $\B(\underline{G},F)$, as defined by Bruhat and
Tits~\cite{bruhat-tits:one,bruhat-tits:two}.  As a $G$-set,
$\B(\underline{G},F)$ is a direct product of an affine space (on which
$G$ acts via translation) and the \emph{reduced building}
$\B^{\text{red}}(\underline{G},F)$, which depends only on
$\underline{G}/\underline{Z}$, where $\underline{Z}$ is the center of
$\underline{G}$.  Note that $Z$ fixes
$\B^{\text{red}}(\underline{G},F)$.  For any extension $E/F$ of finite
residue degree, $\B(\underline{G},F)$ always has a natural embedding
into $\B(\widetilde{\underline{G}},F) = \B(\underline{G},E)$.  To
every point $x\in\B(\underline{G},F)$, there is an associated
parahoric subgroup $G_x$ of $G$.  The stabilizer of $x$ in $G$
contains $G_x$ with finite index.  The pro-$p$-radical of $G_x$ is
denoted $G_{x+}$, and the quotient $G_x/G_{x+}$ is the group of
rational points of a connected reductive $k_F$-group
$\underline{\finiteG}_x$.  These objects depend only on the image of
$x$ in $\B^{\text{red}}(\underline{G},F)$.  Thus, in the case of a
torus $\underline{T}$, we may write $T_0$, $T_{0+}$, and
$\underline{\finiteT}$ instead of $T_x$, $T_{x+}$, and
$\underline{\finiteT}_x$, since these do not depend on the choice of
$x$.  More generally, $G_{0+}$ will denote the set of topologically
unipotent elements in $G$.

We now present an elementary fact about the building that we will use
several times throughout this paper.
\begin{lem}
\label{lem:very_regular}
Let $\underline{Z}$ denote the center of $\underline{G}$, and let
$y,z\in\B(\underline{G},F)$ have distinct images in
$\B^{\textup{red}}(\underline{G},F)$.  Suppose $G_y$ is a maximal
parahoric subgroup, $\gamma\in G_y$, and the image $\bar\gamma$ of
$\gamma$ in $\finiteG_y$ is regular elliptic (i.e., $\bar\gamma$
belongs to no proper $k_F$-parabolic subgroup of 
$\underline{\finiteG}_y$).  Then $\gamma\not\in ZG_z$.
\end{lem}
\begin{proof}
First, suppose $\gamma\in ZG_z\smallsetminus G_z$.
From~\cite[Lemma~4.2.1]{debacker}, $\gamma$ does not fix $z$.
Therefore, $\gamma$ must act on some line containing $z$ via a
nontrivial translation.  By~\cite[Cor.~3.1.5]{debacker}, $\gamma$
cannot fix $y$, a contradiction.

Now suppose $\gamma\in G_z$.  Then $\gamma\in G_x$ for all $x$ lying
on the geodesic between $y$ and $z$.  For such an $x$ that is close
to but not equal to $y$, $G_x$ is a subgroup of $G_y$, and the image
of $G_x$ in $\finiteG_y$ is the group of $k_F$-fixed points of a
proper parabolic subgroup.  Thus $\gamma\not\in G_x$, a contradiction,
and the lemma follows.
\end{proof}

If $\underline{\finiteG}$ is a connected reductive group over a finite
field, $\underline{\finiteT}$ is a maximal torus in $\underline{\finiteG}$,
and $\theta$ is a (complex) character of $\finiteT$, then
let 
$R_\finiteT^\finiteG\theta$
denote the corresponding Deligne-Lusztig virtual character
of $\finiteG$~\cite{deligne-lusztig}.

For any reductive algebraic group $\underline{G}$ defined
over a local or finite field, we have the following notation.
\begin{itemize}
\item
$\mathbf{1}_G$ will denote the trivial representation of $G$.
\item
$\St_G$ will denote the Steinberg representation of $G$.
\item
For any character $\psi$ of $G$, $\St_G(\psi)$ will denote
$\St_G \cdot \psi$.
\item For any representation $\sigma$ of a subgroup $H$ of $G$,
  $\ind_H^G\sigma$ will denote the representation of $G$ obtained from
  $\sigma$ via normalized compact induction.
\item
If $\underline{Z}$ is the center of $\underline{G}$
and $\omega$ is a character of $Z$,
then
let $C(G,\omega)$ denote the space of complex-valued, locally constant
functions $f$ on $G$ such that the support of $f$ is compact modulo $Z$,
and $f(gz)= f(g)\omega(z)$ for all $g\in G$ and $z\in Z$.
\item
$G\reg$ denotes the set of regular semisimple elements of $G$.
\item
For any admissible, finite-length representation
$\pi$ of $G$,
let $\theta_\pi$ denote the character of $\pi$,
considered either as a function on the set of elements or conjugacy classes
of $G$ (of $G\reg$ in the local-field case),
or as a distribution
on an appropriate function space on $G$.
\item Suppose $\varepsilon$ is an automorphism of $G$.  Then
  $\varepsilon$ acts in a natural way on the set of equivalence
  classes of irreducible, admissible representations of $G$.  Suppose
  $\pi$ is such a representation and $\pi\cong \pi^\varepsilon$.  Let
  $\pi(\varepsilon)$ denote an intertwining operator from $\pi$ to
  $\pi^\varepsilon$.  If $\varepsilon$ has order $\ell$, then we can
  and will normalize $\pi(\varepsilon)$ by requiring that the scalar
  $\pi(\varepsilon)^\ell$ equal $1$.  Then $\pi(\varepsilon)$ is well
  determined up to a scalar $\ell$th root of unity.  The
  \emph{$\varepsilon$-twisted character} of $\pi$ is the distribution
  $\theta_{\pi,\varepsilon}$ defined by $\theta_{\pi,\varepsilon}(f) =
  \text{trace}(\pi(f)\pi(\varepsilon))$ for $f\in C_c^\infty(G)$.  As
  with the character, the twisted character can be represented by a
  function (again denoted $\theta_{\pi,\varepsilon}$) on $G$ ($G\reg$
  in the local-field case).  We may regard $\theta_{\pi,\varepsilon}$
  as a function on the set of $\varepsilon$-twisted conjugacy classes.

Note that $\theta_{\pi,\varepsilon}$ still makes sense when $\pi$
is an admissible, finite-length representation.
\item
For any maximal torus $\underline{T}$ of $\underline{G}$,
let $W(T,G)$ denote the quotient of $T$ in its normalizer in $G$,
and let $W_F(\underline{T},\underline{G})$ denote
the group of $F$-points of the absolute Weyl group
$N_{\underline{G}}(\underline{T})/\underline{T}$.
\end{itemize}

\subsection{Shintani lifting}
\label{sec:shintani-general}
Suppose that $E/F$ is a finite, cyclic extension
of local or finite fields,
$\Gamma = \Gal(E/F)$,
and $\underline{G}$ is a connected reductive algebraic $F$-group.
Let $\varepsilon$ denote a generator of $\Gamma$,
and let $\ell$ denote the order of $\Gamma$.
Then one can define a norm mapping from $\widetilde{\underline{G}}$ 
to $\widetilde{\underline{G}}$ by
$$
x \mapsto x \cdot \varepsilon(x) \cdot \,\cdots\, \cdot
\varepsilon^{\ell-1}(x).
$$
If $x$ is defined over $F$ 
then, in general,
the most that one can say about the image of $x$
is that its conjugacy class
in $\underline{G}$ is defined over $F$.
If $F$ is local and $\underline{G}$ has a simply
connected derived group, then such a conjugacy class must
have $F$-points~\cite{rogawski}.
Thus, an $F$-point $x\in\widetilde{G}$
determines a stable conjugacy class in $G$.
Any stable, $\varepsilon$-twisted conjugate of $x$ determines the
same stable conjugacy class in $G$.
Thus, we have a map $\N^{\underline{G}}_{E/F}$
from the set of stable, $\varepsilon$-twisted conjugacy classes
of $\widetilde{G}$ to the set of stable conjugacy classes in $G$.
If $x$ commutes with its Galois conjugates, then we may and will
define $\N^{\underline{G}}_{E/F}(x)\in G$ via the formula above.

Call $g\in\widetilde{G}$ \emph{$\varepsilon$-regular}
if $\N(g)$ is regular.  Let $\widetilde{G}\ereg$ denote
the set of $\varepsilon$-regular elements.

If $\Pi$ and $\widetilde\Pi$ are finite sets of representations
of $G$ and $\widetilde{G}$, respectively, we say that $\widetilde{\Pi}$
is the \emph{Shintani lift} (or \emph{base change}) of 
$\Pi$ if
$$
\Theta_{\widetilde{\Pi},\varepsilon}(g) = 
\Theta_\Pi(\N(g))
$$
for all $g\in \widetilde{G}$ (all $g\in \widetilde{G}\ereg$ in the
local-field case), where $\Theta_\Pi$ and
$\Theta_{\widetilde{\Pi},\varepsilon}$ are nontrivial stable
(resp.~$\varepsilon$-stable) linear combinations of the characters
(resp.~$\varepsilon$-twisted characters) of the elements of $\Pi$
(resp.~$\widetilde{\Pi}$).

If $\underline{T}$ is an $F$-torus
then for any character $\lambda$ of $T$, define the character
$\tilde\lambda$ of $\widetilde T$ by
$\tilde\lambda = \lambda\circ\N^{\underline{T}}_{E/F}$.

\subsection{Notation related to unitary groups}
\label{sec:notation-special}
From now on, fix a nonarchimedean local field $F$ of characteristic
zero with finite residue field $k_F$ of odd order $q$.  Let $E$ be the
unramified quadratic extension of $F$.  Let $E^1$ (resp.~$k_E^1$)
denote the kernel of the norm from $E$ to $F$ (resp.~$k_E$ to $k_F$).
Let $\varepsilon$ denote the nontrivial element of the Galois group
$\Gamma = \Gal(E/F)$.

Let $\underline{G}$ denote a unitary group in three
variables defined with respect to $E/F$.
Then $\underline{G}$
is uniquely determined up to isomorphism, and we can and will assume that
$\underline{G}$ is the unitary group defined by the Hermitian matrix
$$
\Phi =
\left(
\begin{smallmatrix}
0&\phm0&\phm 1\\
0&-1&\phm 0\\
1&\phm0&\phm 0
\end{smallmatrix}
\right).
$$
Then
$$
G = \underline{G}(F) = \left\{ g\in\GL(3,E): g\, \Phi\
  {}^t(\varepsilon(g)) = \Phi\right\}
$$
and $\widetilde{G}=\GL(3,E)$.
Let $\underline{\finiteG}$ denote the corresponding algebraic group over $k_F$.
Then $\widetilde{\finiteG} = \GL(3,k_E)$.

Let $\underline{Z}$ denote the center of $\underline{G}$.
So, following our notational conventions,
$\underline{\widetilde{Z}}$ is the center of $\underline{\widetilde{G}}$.

Let $\B = \B(\underline{G},F)$ and
$\widetilde{\B} = \B(\underline{\widetilde{G}},F) = \B(\underline{G},E)$.
Note that $\varepsilon$ acts on $\widetilde\B$, and we may and will identify
the set of fixed points $\widetilde{\B}^\varepsilon$ with $\B$.

Since $\underline{G}$ is $F$-quasisplit, it contains $F$-Borel
subgroups.
In particular, $\widetilde\B$ must contain some $\varepsilon$-invariant
apartment $\tA$ with more than one $\varepsilon$-fixed point.
Choose an $\varepsilon$-fixed point $y$ in an $\varepsilon$-invariant
minimal facet in $\tA$,
and an $\varepsilon$-invariant alcove $\widetilde\F$ in $\tA$,
such that the closure of $\widetilde\F$ contains $y$.
(Let $\F$ denote the set of $\varepsilon$-fixed points of $\widetilde\F$.)
Then these choices determine an $F$-Borel subgroup $\underline{B}$
together with a Levi factor $\underline{M}$ of $\underline{B}$.
Note that $M$ is isomorphic to $E^\times \times E^1$.
We may assume that our choices of $y$ and $\widetilde\F$ allow us to realize
$\underline{B}$ explicitly as the group of upper triangular matrices
in $\underline{G}$, and $\underline{M}$ as the group
of diagonal matrices.

The boundary of $\F$ contains two points:
the previously chosen point $y$, and another point that we will
denote $z$.
Note that
$\widetilde\F$ is the direct product of a one-dimensional affine space
and an $\varepsilon$-invariant equilateral
triangle $\Delta$ 
in the reduced building of $\widetilde{G}$
(which we will identify with a subset of $\widetilde\B$),
$y$ the $\varepsilon$-fixed vertex of $\Delta$,
and $z$ is the midpoint of the wall of $\Delta$ that is opposite $y$.
In $\B$, $y$ and $z$ are both vertices,
but only $y$ is hyperspecial.

Consider the map $\lambda\colon {U}(1) \rightarrow \underline{G}$
given by $t\mapsto \textup{diag}(1,t,1)$.
Since $\widetilde{\UU(1)} \cong \GL(1)$,
we actually have a one-parameter subgroup of $\widetilde{G}$.
In the usual way, $\lambda$ determines a parabolic $F$-subgroup
$\underline{\widetilde{P}} =\underline{\widetilde{P}}_\lambda$
of $\underline{\widetilde{G}}$, together with a  Levi decomposition
of $\underline{\widetilde{P}}$.
Let $\underline{\widetilde{H}}$ denote the corresponding Levi factor.
Then $\underline{\widetilde{H}}$ is the group of invertible matrices
of the form
$$
\left(
\begin{smallmatrix}
* & 0 & * \\
0 & * & 0 \\
* & 0 & *
\end{smallmatrix}
\right).
$$
This subgroup arises via restriction of scalars from a subgroup
$\underline{H}$ of $\underline{G}$.  Note that $\underline{H}$ is an
$E$-Levi, but not $F$-Levi, subgroup of $\underline{G}$.  It is an
endoscopic group for $\underline{G}$, isomorphic to $\UU(1,1) \times
\UU(1)$.

Similarly, we can define a subgroup $\underline\finiteH$ of
$\underline\finiteG$ and a parabolic $k_F$-subgroup
$\underline{\widetilde\finiteP}$ of $\underline{\widetilde\finiteG}$
with Levi factor $\underline{\widetilde\finiteH}$.  Note that
$\underline\finiteG_y \cong \finiteG$ and
$\underline\finiteG_z \cong \finiteH$.

Up to conjugacy, $\underline{H}$ contains two $F$-tori that are
isomorphic to $\UU(1)\times \UU(1) \times \UU(1)$.  The group of
$F$-points of one of these tori fixes a hyperspecial vertex, and the group
of $F$-points of the other fixes a non-hyperspecial vertex.  Pick such
a torus whose $F$-points fix $y$ (resp.~$z$) and call it $\underline{C}$
(resp.~$\underline{C}'$).  Given the right choices, we can and will
realize $\underline{C}$ as the set of matrices of the form
$$
  \gamma = \left(\begin{array}{ccc}
    \frac{\gamma_1+\gamma_3}{2} & 0 & \frac{\gamma_1-\gamma_3}{2}\\
    0 & \gamma_2 & 0\\
    \frac{\gamma_1-\gamma_3}{2} & 0 & \frac{\gamma_1+\gamma_3}{2}
  \end{array}\right)
$$
where $\gamma_i\in \UU(1)$.  We define the torus
$\underline\finiteC\subset\underline\finiteG$ similarly.  We identify
$\underline{C}$ (and similarly $\underline\finiteC$) with
$\UU(1)\times \UU(1)\times \UU(1)$ via the map $\gamma\mapsto
(\gamma_1 ,\gamma_2 ,\gamma_3 )$.  We will realize $\underline{C}'$ as
$\nu\underline{C}\nu^{-1}$, where
$$
\nu = \left(\begin{array}{ccc}
          1/\sqrt{\smash[b]{\varpi_F}} & 0 & 0\\
          0 & 1 & 0\\
          0 & 0 & \sqrt{\smash[b]{\varpi_F}}
  \end{array}\right).
$$
(The ambiguity in the choice of square root of $\varpi_F$ has
no effect.)

Let $\underline\finiteB_y$ (resp.~$\underline\finiteB_z$) denote the Borel
subgroup of $\underline\finiteG_y$ (resp.~$\underline\finiteG_z$)
determined by $\F$.

For any $F$-group $\underline{L}$, let
$\N^{\underline{L}}=\N^{\underline{L}}_{E/F}$.  When $\underline{L} =
\underline{G}$, we simply write $\N$.  Similarly, for any $k_F$-group
$\underline{\mathsf{L}}$, let $\N^{\underline{\mathsf{L}}} =
\N^{\underline{\mathsf{L}}}_{k_E/k_F}$.  When $\underline{\mathsf{L}} =
\underline{\mathsf{G}}$, we simply write $\bar\N$.

For any subgroup $\underline{S}\subset\underline{\widetilde{G}}$, let
$\det_{\underline S}$ denote the restriction of the determinant
to $\underline{S}$.  We will omit the subscript when it is
clear from the context.  Similar notation holds for subgroups
of $\underline{\widetilde{\finiteG}}$.

\subsection{Cartan subgroups of $\underline{G}$}
\label{sec:cartan-G}
For a quadratic extension $L/K$, denote by $\UU(1,L/K)$ the
unitary group in one variable over $K$ defined with respect to $L/K$.
Up to stable conjugacy, there are four kinds of Cartan subgroup of
$\underline{G}$.  In the notation of~\cite{rogawski}, they are
isomorphic to:
\begin{enumerate}[(\ref{sec:cartan-G}--1)]
\addtocounter{enumi}{-1}
\item $R_{E/F}(\GL(1))\times \UU(1,E/F)$,
\item $\UU(1,E/F)\times \UU(1,E/F)\times \UU(1,E/F)$,
\item $R_{E/F}(\UU(1,EK/K))\times \UU(1,E/F)$ for $K$ a ramified quadratic
  extension of $F$,
\item $R_{L/F}(\UU(1,EL/L))$ for $L$ a cubic extension of $F$.
\end{enumerate}

\subsection{Representation theory of $\finiteG$, $\finiteH$, and $\finiteC$}
\label{sec:reps-finite}
A reference for much of this section is~\cite{srinivasan}.

\subsubsection*{Representations of $\finiteG$}
Let $\finiteB$ denote a Borel subgroup of $\finiteG$
with Levi factor $\finiteM$,
and let $\theta$ be a character of $\finiteM$.
Then the induced representation $\ind_\finiteB^\finiteG \theta$
is irreducible except when $\theta$ extends to a character
of $\finiteH$.  
In this case, the induced representation is a sum of two irreducible
components.  If $\theta$ extends to a character $\theta_0$ of $\finiteG$, then 
these components are $\theta_0$ and $\St_\finiteG(\theta_0)$.

Let $L$ denote a cubic unramified extension of $E$.  Then $\finiteG$
contains a torus $\finiteS$ that is isomorphic to the kernel of the
norm map from $k_{EL}$ to $k_L$.  Let $\finiteT$ be either $\finiteS$
or $\finiteC$.  For any character $\theta$ of $\finiteT$ with trivial
stabilizer in $W_{k_F}(\underline{\finiteT},\underline{\finiteG})$, we
have a Deligne-Lusztig cuspidal representation whose character is
$-R_\finiteT^\finiteG \theta$.  For $\finiteT = \finiteS$, we will
call such representations ``cubic cuspidal representations.''

The other irreducible representations of $\finiteG$ have the form
$\tau\cdot\psi$, where $\tau$ is the cuspidal unipotent representation
and $\psi$ is a character.

\subsubsection*{Representations of $\finiteH$}
As above let $\finiteB$ be a Borel subgroup of $\finiteG$ with Levi
factor $\finiteM$ and let $\theta$ be a character of $\finiteM$.
The induced representation $\ind_{\finiteB\cap \finiteH}^\finiteH \theta$
is irreducible except when $\theta$ extends to a character
$\theta_1$ of $\finiteH$.  In this case, the induced representation is
the sum of $\theta_1$ and $\St_\finiteH(\theta_1)$.

The remaining representations of $\finiteH$ are the Deligne-Lusztig
cuspidal representations, whose characters are of the form
$-R_\finiteC^\finiteG \theta$ for $\theta\in\Hom (\finiteC ,\C^\times)$ in
general position with respect to the action of
$W_{k_F}(\underline{\finiteC},\underline{\finiteH})$.

\subsubsection*{Representations of $\finiteC$}
We will need a technical result on linear combinations of characters
of $\finiteC$.  Let $A$ be a finite abelian group of order $n$, and
let $\chi_1 ,\ldots ,\chi_n$ be the irreducible characters of $A$.
The following three lemmas concern characters of products of copies of
$A$.

\begin{lem}
\label{lem:cyclic}
  Let
  $$
  f=\sum_{i=1}^n a_i \chi_i,
  $$
  where $a_i\in\C$.  Suppose that $f$
  vanishes off of a subset of $A$ of size $2$.  Then either
  $f=0$, or the number of $i$ such that $a_i \neq 0$ is at least
  $n/2$.
\end{lem}
\begin{proof}
  Let $\{ a,b\}$ be the above subset of $A$.  We have
$$
na_i = n\cdot \langle f,\chi_i\rangle
= f(a)\bar{\chi}_i(a) + f(b)\bar{\chi}_i(b).
$$
  Assume $f\neq 0$.  If $f(b)=0$, then for all $i$, $a_i =
  f(a)\bar{\chi}_i(a)/n\neq 0$.
  If $f(b)\neq 0$, then $a_i\neq0$ unless $\bar{\chi}_i(ba^{-1}) =
  -f(a)/f(b)$.  Since $ba^{-1}\neq 1$, this equality holds for at most
  $n/2$ values of $i$.
\end{proof}

\begin{lem}
\label{lem:cyclic2}
  Let $N$ be the subset of
  $A\times A$ consisting of all elements $(a,b)$ such that
  $a\neq b$.  Suppose that for some $a_{ij}\in\C$,
  $$
  f=\sum_{i,j} a_{ij} \chi_i\otimes\chi_j
  $$
  vanishes on $N$.  Then either $f=0$ or at
  least $n$ of the $a_{ij}$ are nonzero.
\end{lem}
\begin{proof}
  Assume $f\neq 0$.  Fix $a\in A$.  Evaluating $f$ at $(a,b)$ for
  $b\neq a$, we obtain that the function
  $$
  \sum_j\left( \sum_i a_{ij}\chi_i(a)\right) \chi_j
  $$
  on $A$
  vanishes on $A-\{ a\}$.  It follows easily that either this function
  vanishes on $A$, or for all $j$, the coefficient $\sum_i
  a_{ij}\chi_i(a)$ is nonzero.  The former case cannot happen since
  $f\neq 0$.  In the latter case, it follows that for all $j$, at
  least one coefficient $a_{ij}$ must be nonzero.  Hence at least $n$
  of the $a_{ij}$ must be nonzero.
\end{proof}

\begin{lem}
  Let $N'$ be the subset of $A\times A\times A$ consisting of all
  elements $(a,b,c)$ such that $a$, $b$, and $c$ are distinct.
  Suppose that for some $a_{ijk}\in\C$,
  $$
  f= \sum_{i,j,k}a_{ijk}\chi_i\otimes\chi_j\otimes\chi_k
  $$
  vanishes on $N'$.  Then either $f$ vanishes on $A\times A \times A$
  or at least $n/2$ of the $a_{ijk}$ are nonzero.
\end{lem}

\begin{proof}
Assume $f\neq 0$.
Fix $a\neq b$ in $A$.  Then the function
$$
\sum_k\biggl(\sum_{i,j}a_{ijk}\chi_i(a)\chi_j (b)\biggr)\chi_k
$$
on $A$ vanishes off of $\{ a,b \}$.  Hence, by
Lemma~\ref{lem:cyclic}, either this function vanishes on $A$, or for
at least $n/2$ values of $k$, the coefficient
$\sum_{i,j}a_{ijk}\chi_i(a)\chi_j (b)$ is nonzero.  In the latter
case, for each such $k$, at least one coefficient $a_{ijk}$ must be
nonzero.  Hence at least $n/2$ of the $a_{ijk}$ are nonzero.

We may therefore assume that the former case holds for all pairs
$a\neq b$.  By the linear independence of characters, the coefficient
$\sum_{i,j}a_{ijk}\chi_i(a)\chi_j (b)$ must vanish for all $k$ and all
pairs $a\neq b$.  Since $f\neq 0$, $a_{i'j'k'}\neq 0$ for some
$i',j',k'$.  Thus the function $\sum_{i,j}a_{ijk'}\chi_i\otimes \chi_j
$ on $A\times A$ vanishes on the set $N$ of Lemma~\ref{lem:cyclic2},
but it does not vanish on $A\times A$ since $a_{i'j'k'}\neq 0$.  Hence
Lemma~\ref{lem:cyclic2} implies that at least $n$ of the coefficients
$a_{ijk'}$ must be nonzero.
\end{proof}

\begin{cor}
\label{cor:lin_dep}
Suppose that 
$$
\sum_{\mbox{$\chi\in\Hom (\finiteC ,\C^\times )$}}a_{\chi}\chi
$$
vanishes on $\finiteC \cap \finiteG\reg$, where
$a_{\chi}\in\C$.  Then either this linear combination vanishes on
$\finiteC$ or at least $(q+1)/2$ of the $a_{\chi}$ are nonzero.
\qed
\end{cor}

\begin{cor}
\label{cor:q-value}
Suppose that $q>59$, and
let $f = \sum a_{\chi}\chi$
be a linear combination of at most 30 characters of $\finiteC$.
If $f$ vanishes on $\finiteC\cap \finiteG\reg$,
then $f$ vanishes on $\finiteC$.
\qed
\end{cor}

\subsection{Shintani lifting for $\finiteG$ and $\finiteH$}
\label{sec:Shintani-finite}
According to~\cite{srinivasan}, the irreducible characters of
$\finiteG$ are of the form $\pm R_\finiteL^\finiteG\theta$, where
$\underline\finiteL$ is the connected centralizer of some semisimple
element of $\underline\finiteG$, and $\theta$ is the twist of a
unipotent character of $\finiteL$ by a one-dimensional character in
general position.  Moreover, one obtains a cuspidal character of
$\finiteG$ precisely when $\underline\finiteL$ is an elliptic torus or
when $\underline\finiteL=\underline\finiteG$ and $\theta$ is a twist
of the unique cuspidal unipotent character of $\finiteG$.

By~\cite{kawanaka}, our assumption that $k_F$ has odd characteristic
guarantees the existence of Shintani descent from $\widetilde\finiteG$
to $\finiteG$.  In~\cite{digne}, Digne gives a general proof that
Shintani descent is compatible with Deligne-Lusztig induction.  In
particular, if $\sigma$ is an irreducible representation of $\finiteG$
with character $\pm R_\finiteL^\finiteG\theta $ ($\theta$ a character
of $L$), then the character of the Shintani lift $\tilde\sigma$ of
$\sigma $ from $\finiteG$ to $\widetilde\finiteG$ is of the form $\pm
R_{\widetilde\finiteL}^{\widetilde\finiteG}\tilde{\theta}$, where
$\tilde{\theta}$ is the Shintani lift of $\theta$.  Now
$\underline{\widetilde\finiteL}$ is a Levi factor of a parabolic
subgroup of $\underline{\widetilde\finiteG}$ unless $\finiteL$ is
isomorphic to the torus $\finiteS$ defined in \S\ref{sec:reps-finite}.
Hence $\tilde\sigma$ is a parabolically induced representation unless
$\finiteL\cong\finiteS$ or $\finiteL =\finiteG$.  In the former case,
$\widetilde\finiteS$ is an elliptic torus isomorphic to
$k_{EL}^\times$ and $\tilde\sigma$ is cuspidal.  In the latter case,
$\sigma $ is a one-dimensional representation
$\varphi\circ\det_{\underline{\finiteG}}$, a twist $\St_\finiteG
(\varphi\circ\det_{\underline{\finiteG}})$ of the Steinberg
representation, or a twist
$\tau(\varphi\circ\det_{\underline{\finiteG}})$ of the cuspidal
unipotent representation.  One shows easily that the Shintani lift
$\tilde\sigma$ is, respectively,
$\tilde\varphi\circ\det_{\underline{\widetilde\finiteG}}$,
$\St_{\widetilde\finiteG}(\tilde\varphi\circ\det_{\underline{\widetilde\finiteG}})$,
or $\tilde\tau
(\tilde\varphi\circ\det_{\underline{\widetilde\finiteG}})$, where
$\tilde\tau$ is the unipotent representation of $\widetilde\finiteG$
not equivalent to ${\bf 1}_{\widetilde\finiteG}$ or
$\St_{\widetilde\finiteG}$.  The remaining representations of
$\finiteG$ are those whose characters are of the form
$R_\finiteH^\finiteG\theta$.  By~\cite{digne}, the Shintani lifts of
such representations are representations induced from
$\widetilde\finiteP$.  Hence the cubic cuspidal representations of
$\finiteG$ are exactly those irreducible representations of $\finiteG$
whose Shintani lifts are cuspidal.

We now consider Shintani lifting for irreducible representations
of $\finiteH$.
From \S\ref{sec:reps-finite},
most such representations have characters of the form 
$\pm R_\finiteT^\finiteH \theta$.
From Digne~\cite{digne},
the Shintani lift of such a representation has character
$\pm R_{\widetilde{\finiteT}}^{\widetilde{\finiteH}} \widetilde\theta$.

The remaining representations of $\finiteH$ are the one-dimensional
representations $\varphi\circ \det_{\underline{\finiteH}}$ and the Steinberg
representations $\St_\finiteH(\varphi\circ\det_{\underline{\finiteH}})$.
It is easy to see that the respective Shintani lifts of these
representations are $\widetilde\varphi\circ\det_{\underline{\widetilde{\finiteH}}}$
and $\St_{\finiteH}(\widetilde\varphi\circ\det_{\underline{\widetilde{\finiteH}}})$.

\subsection{Depth-zero representations of $G$}
\label{sec:depth-zero-G}

\subsubsection*{Principal series of $G$}

For $\lambda\in\Hom (M,\C^\times)$, there exist
unique characters $\lambda_1\in\Hom (E^{\times},\C^\times)$ and
$\lambda_2\in\Hom (E^1,\C^\times)$ such that 
\begin{equation}
\label{eqn:lambda}
\lambda\left(\left(
\begin{smallmatrix}
          \alpha & 0 & 0\phinv \\
          0 & \beta & 0\phinv \\
          0 & 0 & \bar{\alpha}^{-1}
  \end{smallmatrix}
\right)\right)  =
\lambda_1(\alpha)\lambda_2(\alpha\bar{\alpha}^{-1}\beta ),
\end{equation}
where $\alpha\in E^\times$, $\beta\in E^1$.
By~\cite{Keys}, $\ind_B^G\lambda$ is
irreducible except for in the following cases:
\begin{enumerate}[(\ref{sec:depth-zero-G}PS--1)]
\item
$\lambda_1 = |\cdot|_E^{\pm 1}$
\item
$\lambda_1|_{F^\times} = \omega_{E/F} \, |\cdot |_F^{\pm 1}$
\item
$\lambda_1$ is nontrivial and $\lambda_1|_{F^\times}$ is trivial.
\end{enumerate}

In case~\reducible{1},
$\ind_B^G\lambda$ has two constituents: the
one-dimensional representation $\psi=\lambda_2\circ\det$, and the
square-integrable Steinberg representation $\St_G (\psi)$.

In case~\reducible{2},
$\ind_B^G\lambda$ also has two constituents: a square-integrable
representation $\pi^2(\lambda)$ and a non-tempered unitary
representation $\pi^n(\lambda)$.

In case~\reducible{3},
$\ind_B^G\lambda$ decomposes into a direct sum
$\pi_1(\lambda)\oplus\pi_2 (\lambda)$.

By~\cite{moy-prasad2}, $\ind_B^G\lambda$ has depth zero if and only if
$\lambda$ has depth zero.

\subsubsection*{Other representations of $G$}
Since $G$ has no non-minimal proper parabolic subgroups, the remaining
irreducible representations are all supercuspidal.  From
either~\cite{moy-prasad2} or~\cite{morris:intertwining}, we know that
all such representations have a unique expression of the form
$\ind_{G_x}^G \sigma$, where $x=y$ or $z$, and $\sigma$ is the
inflation to $G_x$ of an irreducible cuspidal representation
$\overline\sigma$ of $\finiteG_x$.  The representations
$\overline\sigma$ are classified in \S\ref{sec:reps-finite}.  Based on
this classification, we have the following kinds of supercuspidal
representation of depth zero.
\begin{enumerate}[(\ref{sec:depth-zero-G}SC--1)]
\item $\ind_{G_y}^G\sigma$, where $\bar\sigma$ is a cubic cuspidal
  representation of $\finiteG_y\cong \finiteG$.
\item $\ind_{G_y}^G\sigma$, where $\bar\sigma$ is a cuspidal
  representation of $\finiteG_y$ with character
$-R_{\finiteC}^{\finiteG_y}\varphi$ and
  $\varphi=\varphi_1\otimes\varphi_2\otimes\varphi_3$ is a regular character of
  $\finiteC$
(with respect to $W_{k_F}(\underline{\finiteC} ,\underline{\finiteG}_y)$).
\item $\ind_{G_y}^G\sigma$, where $\bar\sigma$ is the twist $\tau\cdot
  (\eta\circ\det)$ of the cuspidal unipotent representation
$\tau$ of $\finiteG_y$,
  and $\eta\in\Hom (k_E^1,\C^\times)$.
\item $\ind_{G_z}^G\sigma$, where $\bar\sigma$ is a cuspidal
  representation of $\finiteG_z\cong\finiteH$ with character
  $-R_{\finiteC'}^{\finiteG_z}\varphi$ and
  $\varphi=\varphi_1\otimes\varphi_2\otimes\varphi_3$ is a regular character of
  $\finiteC'$
(with respect to $W_{k_F}(\underline{\finiteC}' ,\underline{\finiteG}_z)$).
Recall that, according to our notational conventions,
$\underline{\finiteC}'$
is the finite reductive quotient of the (unique) parahoric
subgroup of $\underline{C}'$.
\end{enumerate}

\section{Description of depth-zero $L$-packets and explicit base
  change for unitary groups in two variables}
\label{sec:u11}

In this section we give brief descriptions of the depth-zero
$L$-packets for the quasi-split
group $U(1,1)(F)$ and the compact group $U(2)(F)$,
as well as their base change lifts to $\GL_2(E)$.
We omit the proofs as they are entirely analogous to
(but less complicated than) those for $U(2,1)$.  Let $\underline{H}^0$
be the group $U(1,1)$, which we will view as the subgroup of
$\underline H$
consisting of all matrices of the form
$$
\left(
\begin{smallmatrix}
* & 0 & * \\
0 & 1 & 0 \\
* & 0 & *
\end{smallmatrix}
\right).
$$
For every subgroup $\underline{L}$ of $\underline{G}$, let
$\underline{L}^0$ denote the subgroup
$\underline{L}\cap\underline{H}^0$ of $\underline{H}^0$.
Let $\underline{H}^1$ denote the compact inner form $U(2)$ of $H^0$.
Note that
$\underline{\widetilde{H}}^0(F)
\cong \underline{\widetilde{H}}^1(F) \cong \GL_2(E)$.

From our descriptions, it will be clear that the analogue
of Theorem~\ref{thm:main} holds for unitary groups
in two variables.

\subsection{Depth-zero $L$-packets for $U(1,1)$}
\label{sub:L-packets_u11}
The $L$-packets of $H^0$ are the ${\rm PGL}_2(F)$-orbits on the
set of equivalence classes of irreducible admissible representations
of $H^0$~\cite[\S11.1]{rogawski}.  We first describe the
principal series $L$-packets.

Let $\lambda\in\Hom (M^0,\C^\times) = \Hom (E^\times ,\C^\times )$.
According to~\cite[\S11.1]{rogawski}, the principal series
$\ind_{B^0}^{H^0}\lambda$ is irreducible except in the cases
\begin{enumerate}
\item
$\lambda |_{F^\times} = |\cdot|_F^{\pm 1}$
\item
$\lambda |_{F^\times} = \omega_{E/F} $.
\end{enumerate}

In the first case, $\ind_{B^0}^{H^0}\lambda$ has two constituents: the
one-dimensional representation $\psi = \mu\circ\det$, where
$\mu\circ\N = \lambda |\cdot |_E^{\mp 1/2}$, and the Steinberg
representation $\St_G (\psi)$.  In the second case,
$\ind_{B^0}^{H^0}\lambda$ decomposes into a direct sum
$\pi_1(\lambda)\oplus\pi_2 (\lambda)$
of irreducible representations.
By~\cite{moy-prasad2},
$\ind_{B^0}^{H^0}\lambda$ has depth zero if and only if $\lambda$ has
depth zero.

The principal series $L$-packets of $G$ are as follows~\cite[\S
11.1]{rogawski}.  (Here $\lambda$ and $\psi$ denote one-dimensional
representations of $M^0$ and $H^0$, respectively.)
\begin{enumerate}
\item $\{ \ind_{B^0}^{H^0}\lambda\}$, where $\ind_{B^0}^{H^0}\lambda$
  is irreducible;
\item $\{\psi\}$; 
\item $\{\St_{H^0} (\psi)\}$;
\item $\{\pi_1 (\lambda),\pi_2 (\lambda)\}$, where
  $\ind_{B^0}^{H^0}\lambda$ is reducible of the second type described above.
\end{enumerate}

The remaining irreducible representations and $L$-packets of $G$ are
all supercuspidal.  The depth zero supercuspidals of $H^0$ have a
unique expression of the form $\ind_{H^0_v}^{H^0} \sigma$, where $v=y$
or $z$, and $\sigma$ is the inflation to $H^0_v$ of an irreducible
cuspidal representation $\overline\sigma$ of $\finiteH^0_v$.  Let
$\underline{C}^0_v$ be $\underline{C}\cap \underline{H}^0$ if $v=y$,
and $\underline{C}'\cap \underline{H}^0$ if $v=z$.  Then the character
of such a representation $\overline\sigma$ must be of the form
$-R_{\finiteC^0_v}^{\finiteH^0_v}\varphi$, where $\varphi$ is a
character of $\finiteC^0_v$ in general position.
Since $\finiteC^0_v \cong k_E^1\times k_E^1$, we may
view any such character as having the form $\varphi_1\otimes\varphi_2$,
where the $\varphi_i$ are distinct characters of $k_E^1$.

Fix a cuspidal representation $\overline\sigma$ of $\finiteH^0$.
Viewing it as a representation of $\finiteH^0_v$, we inflate it to a
representation $\sigma_v$ of $H^0_v$.  Let $\pi_v =
\ind_{H^0_v}^{H^0}\sigma_v$.  Then $\{\pi_y ,\pi_z\}$ is a depth-zero
supercuspidal $L$-packet of $H^0$.  Conversely, all such $L$-packets
are of this form.
If $\pi_y$ and $\pi_z$ are formed from the character
$\varphi_1\otimes \varphi_2$ of $\finiteC^0$ as above,
then for future reference call this $L$-packet
$\Pi^0_{\varphi_1,\varphi_2}$.

\subsection{Base change lifts for $U(1,1)$}
\label{sub:bc_u11}
By~\cite[\S11.4]{rogawski}, the base change lifts of principal series
$L$-packets of $H^0$ are as follows.  Let $\lambda\in\Hom
(M^0,\C^\times)$.
\begin{enumerate}[(i)]
\item
If $\ind_{B^0}^{H^0}\lambda$ is irreducible and
  $\ind_{\widetilde B^0}^{\widetilde H^0}\tilde\lambda$ is irreducible,
  then the base
  change lift of the $L$-packet $\{ \ind_{B^0}^{H^0} \lambda\}$ is
  $\ind_{\widetilde{B}^0}^{\widetilde{H}^0}\tilde\lambda$.
\item
If $\ind_{B^0}^{H^0}\lambda$ is irreducible but
  $\ind_{\widetilde B^0}^{\widetilde H^0}\tilde\lambda$ is reducible,
  then $\lambda|_{F^\times} = |\phantom{x}|_F^{\pm 1}\omega_{E/F}$, and
  the base
  change lift of the $L$-packet $\{ \ind_{B^0}^{H^0} \lambda\}$ is
  $\lambda |\phantom{x}|_E^{\mp 1/2} \circ \det$.
\item If $\lambda|_{F^\times} = |\cdot |_F^{\pm 1}$, let $\psi$ be the
  one-dimensional representation $\mu\circ\det_{\underline{H}^0}$,
  where $\mu\circ\N = \lambda |\cdot |_E^{\mp 1/2}$.  Then the lift of
  the $L$-packet consisting of the constituent $\psi$ (resp., the
  Steinberg constituent $\St_{H^0}(\psi)$) of
  $\ind_{B^0}^{H^0}\lambda$ is the one-dimensional constituent
  $\tilde{\psi} = (\lambda |\cdot |_E^{\mp
    1/2})\circ\det_{\underline{\widetilde{H}}^0}$ (resp., the
  Steinberg constituent $\St_{\widetilde{H}^0}(\tilde{\psi})$) of
  $\ind_{\widetilde{B}^0}^{\widetilde{H}^0}\tilde\lambda$.
\item If $\lambda|_{F^\times} = \omega_{E/F}$, then the lift of the
  $L$-packet $\{\pi_1(\lambda),\pi_2(\lambda)\}$ is
  $\ind_{\widetilde{B}^0}^{\widetilde{H}^0}\tilde\lambda$.
\end{enumerate}

The base change lift of
the depth-zero supercuspidal $L$-packet
$\Pi^0_{\varphi_1,\varphi_2}$ 
is the principal series representation
$\ind_{\widetilde{B}^0}^{\widetilde{H}^0}\varphi^*$,
where $\varphi^*$
is the character $\hat\varphi_1\omega_{E'/E}\otimes\hat\varphi_2\omega_{E'/E}$
of $E^\times\times E^\times\cong\widetilde{M}^0$.  Here, $E'$ is an
unramified quadratic extension of $E$,
and $\hat\varphi_i$ is the inflation to $E^\times$ of the character
$\tilde\varphi_i$ of $k_E^\times$.

\subsection{Depth-zero $L$-packets for $U(2)$}
Since $H^1$ is compact, it has only one parahoric subgroup
(and in fact is equal to it).
The finite reductive quotient $\finiteH^1$ is isomorphic to
$k_E^1 \times k_E^1$.  Thus, every irreducible,
depth-zero representation of $H^1$ has the form
$\infl(\varphi_1\otimes\varphi_2)$,
the inflation to $H^1$ of a character of $\finiteH^1$.

Let
$$
\Pi^1_{\varphi_1,\varphi_2} =
\begin{cases}
\{\infl(\varphi_1\otimes\varphi_2),
\infl(\varphi_2\otimes\varphi_1)\}
& \text{if $\varphi_1 \neq \varphi_2$}, \\
\{\infl(\varphi_1\otimes\varphi_2)\}
& \text{if $\varphi_1 = \varphi_2$}.
\end{cases}
$$
Then we declare the $\Pi^1_{\varphi_1,\varphi_2}$
to be the $L$-packets for $H^1$.
These $L$-packets are chosen so as to make the
correspondence $JL$
given in \S\ref{sec:bc-u2} work properly.

\subsection{Base change lifts for $U(2)$ via a Jacquet-Langlands-like
correspondence}
\label{sec:bc-u2}
Since $\underline{H}^1$ is an inner form of $\underline{H}^0$,
we can obtain a base
change lift if we can associate each $L$-packet for $H^1$
to one for $H^0$.  This association will be similar to
the Jacquet-Langlands correspondence (or ``Abstract Matching
Theorem''~\cite{dkv:amt,rogawski:amt,badulescu:thesis}).
That is, given an $L$-packet $\Pi^1$ for $H^1$, we want to 
find an $L$-packet $\Pi^0$ for $H^0$ such that
\begin{equation}
\label{eqn:jl-u2}
\sum_{\pi\in\Pi^1} \theta_\pi(g_1)
=
\pm
\sum_{\pi\in\Pi^0} \theta_\pi(g_0)
\end{equation}
for all regular $g_1\in H^1$
and $g_0\in H^0$ whose stable conjugacy classes are associated
in a natural way.

Define a map $JL$ from the depth-zero $L$-packets of $H^1$
to those of $H^0$ by
$$
JL(\Pi^1_{\varphi_1,\varphi_2}) =
\Pi^0_{\varphi_1,\varphi_2}
$$
if $\varphi_1\neq\varphi_2$.
If $\varphi = \varphi_1=\varphi_2$, then we define
$JL(\Pi^1_{\varphi_1,\varphi_2})$
as follows.
Form the character $\varphi\circ\N$ of $k_E^\times$,
which we can then inflate to a character $\lambda$ of $E^\times$.
Now let
$JL(\Pi^1_{\varphi_1,\varphi_2})$
be the Steinberg component of
$\ind_{B^0}^{H^0} \lambda  |\cdot |_F$.
More specifically, this representation is
$\St_{H^0}(\mu\circ\det)$,
where $\mu\circ\N = \lambda |\cdot|_F$,
as in \S\ref{sub:L-packets_u11}.

It is not difficult to see that $JL$ is the only correspondence
that satisfies
\eqref{eqn:jl-u2}
for all $g_1\in H^1$
whose image in $\finiteH^1$ is regular.
Thus, if we assume that there is a Jacquet-Langlands-like
correspondence from the depth-zero $L$-packets of $H^0$
to those of $H^1$, then it must be $JL$.
\section{Description of depth-zero $L$-packets and $A$-packets for $G$}
\label{sec:L-packets}

In almost all cases, $L$-packets and $A$-packets are the same.
In one case (see below), a certain principal series $L$-packet is enlarged
to form an $A$-packet.  Thus, while the $L$-packets constitute
a partition of the set of equivalence classes of irreducible representations,
the $A$-packets do not.

\subsection{$L$-packets consisting of principal series constituents}
\label{sub:ps}

The following proposition is due to Rogawski~\cite[\S12.2]{rogawski}.

\begin{prop}
\label{prop:ps_L-packets}
The $L$-packets of $G$ that consist entirely of principal series
constituents all have one of the following forms
(where $\lambda$ and $\psi$ denote one-dimensional representations
of $M$ and $G$, respectively):
\begin{enumerate}[(\ref{sub:ps}--1)]
\item  $\{ \ind_G^P\lambda\}$,
  where $\ind_G^P\lambda$ is irreducible;
\item $\{\psi\}$; 
\item
$\{\St_G (\psi)\}$;
\item $\{\pi_1 (\lambda),\pi_2 (\lambda)\}$, where $\ind_G^P\lambda$ is
  reducible of type \textup{\reducible{3}}.
\item $\{\pi^n (\lambda)\}$, where $\ind_G^P\lambda$ is
  reducible of type \textup{\reducible{2}}.
\end{enumerate}
\end{prop}

In the last case, $\pi^n (\lambda)$
is contained in the $A$-packet $\Pi (\lambda) =
\{ \pi^n (\lambda),\pi^s (\lambda)\}$, where $\pi^s (\lambda)$
is the supercuspidal
representation that sits inside an $L$-packet with the 
square-integrable principal series constituent $\pi^2(\lambda)$.  In
the depth-zero setting, the representation $\pi^s (\lambda)$ will be
explicitly described in~\S\ref{sub:non-stable}

\subsection{Singleton supercuspidal $L$-packets}
\label{sub:stable_sc}

In this section, we characterize the stable supercuspidal
representations of $G$ of depth zero in terms of inducing data.  

\begin{prop}
\label{prop:stable_sc}
  A supercuspidal representation $\pi$ of $G$ of depth zero is stable
  if and only if $\pi$ is of the form $\ind_{G_y}^G\sigma$, where
  $\sigma$ is the inflation to $G_y$ of a cubic cuspidal
  representation $\bar\sigma$ of $\finiteG_y$.
\end{prop}
\begin{proof}
Let $\pi$ be a representation of the above form.  Let $\gamma$ be an
element of $G\reg$ and let $\gamma'$
be a stable conjugate of $\gamma$.
We will show that
\begin{equation}
\label{eq:stable}
\theta_\pi (\gamma)=\theta_\pi (\gamma').
\end{equation} 
The conjugacy classes contained within the stable conjugacy class
of $\gamma$ are parametrized by
$$
\text{Ker}\{H^1(F,\underline{G}_\gamma)\rightarrow H^1(F,\underline{G})\}
$$
(see~\cite[\S3.1]{rogawski}),
where $\underline{G}_\gamma$ is the centralizer of
$\gamma$ in $\underline{G}$.
If $\gamma$ is contained in a Cartan subgroup of $G$
of type \cartan{0} or \cartan{3},
then this kernel is trivial by~\cite[\S3.6]{rogawski}
so any stable
conjugate $\gamma'$ of $\gamma$ is a conjugate of $\gamma$.  Hence
$\theta_\pi (\gamma)=\theta_\pi (\gamma')$.  Therefore, we may assume
that $\gamma$ is contained in a Cartan subgroup $T$ of type
\cartan{1} or \cartan{2}.  

For any regular, depth-zero $X$ in the dual of the Lie algebra of a
cubic torus in $G$, the germ $\theta_{\pi}|_{G_{0+}}$ coincides with a
constant multiple of the Fourier transform of the orbital integral
corresponding to $X$.  (This follows from Corollaire III.10 and
Proposition III.8 of~\cite{waldspurger:nilpotent}.  It also follows
from the proof of the main theorem
of~\cite{adler-debacker:mk-theory}.)  The Weyl group of a cubic torus
acts via the Galois group, so two regular elements of the torus are
conjugate if and only if they are stably conjugate.  Moreover, every
stable conjugate of a cubic torus is conjugate to it.  Therefore, this
orbital integral is stable.  From~\cite{waldspurger:fourier}, the
Fourier transform of a stable distribution is stable.  Thus,
$\theta_\pi|_{G_0+}$ is stable.  If $z\in Z$, it is clear that
$\theta_\pi (\gamma z) = \theta_\pi (\gamma' z)$ if and only if
$\theta_\pi (\gamma) = \theta_\pi (\gamma')$.  Thus,
$\theta_\pi|_{ZG_0+}$ is stable.

It follows that (\ref{eq:stable}) holds if $\gamma\in ZT_{0+}$.
Therefore, suppose that $\gamma\notin ZT_{0+}$.
We will show that $\theta_\pi$ vanishes at all stable conjugates
of $\gamma$ (including $\gamma$ itself), thus establishing
\eqref{eq:stable}.
Let $\gamma''$ be a
stable conjugate of $\gamma$.
If no conjugate of $\gamma''$ is contained in $G_y$,
then $\theta_\pi(\gamma'') = 0$ from Proposition~\ref{induced_char}.
So assume $\gamma''\in G_y$. It follows
easily from our assumptions on $\gamma$ that the
characteristic polynomial of the image $\bar{\gamma}''$ of $\gamma''$
in $\finiteG_y$ is reducible over $k_E$ and that its roots are not all
the same.  But then the semisimple part of $\bar{\gamma}''$ is not
contained in a cubic torus of $\finiteG_y$ so,
by~\cite[6.9]{srinivasan}, it follows that
$\theta_{\bar\sigma}(\bar{\gamma}'')=0$.
Thus $\theta_\pi (\gamma'') = 0$ by Proposition~\ref{induced_char}.

Conversely, suppose that $\pi$ is not of the form given
in the statement of the proposition.  By the
classification in \S\ref{sec:depth-zero-G}, it follows that $\pi$
is of type~\depthzerosc{2},
\depthzerosc{3},
or
\depthzerosc{4}.

Let $\gamma=(\gamma_1, \gamma_2, \gamma_3)\in C\subset G_y$ have
regular image $\bar{\gamma}$ in $\finiteG_y$.  Let $\gamma'\in G_z$ be
the conjugate of $\gamma$ by $\nu$ (see~\S\ref{sec:notation-special}).
Then $\gamma\in G\reg$, and $\gamma$ lies in a unique maximal
parahoric by Lemma~\ref{lem:very_regular}, namely $G_y$.  Also, the image
$\bar\gamma'$ of $\gamma'$ in $\finiteG_z$ is regular elliptic, so
that $\gamma'$ is not contained in any parahoric other than $G_z$.  We
note that $\gamma$ and $\gamma'$ are stably conjugate elements of $G$
that are not conjugate in $G$.

If $\pi$ is of type \depthzerosc{2} or \depthzerosc{3},
then $\pi$ is compactly induced
from the inflation $\sigma$ to $G_y$ of a non-cubic cuspidal
representation $\bar\sigma$ of $\finiteG_y$.
Thus $\theta_\pi (\gamma')=0$ by
Proposition~\ref{induced_char} since $\gamma'$ is not contained in any conjugate
of $G_y$.  On the other hand, since the only conjugate of $G_y$
containing $\gamma$ is $G_y$,
$
\theta_{\pi} (\gamma) = \theta_\sigma (\bar{\gamma})
$
by Proposition~\ref{induced_char}.

Suppose that $\pi$ is of type \depthzerosc{3},
i.e., $\bar\sigma=\tau\cdot
    (\eta\circ\det)$ where $\tau$ is the cuspidal unipotent
    representation of $\finiteG_y$ and $\eta\in\Hom (k_E^1,\C^\times)$.  Then
$$
\theta_{\bar\sigma} (\bar{\gamma}) = 2(\eta\otimes\eta\otimes\eta)
  (\bar{\gamma})\neq 0
$$
  by~\cite[p.~31]{ennola}.  On the other hand, suppose that $\pi$
  is of type \depthzerosc{2}, i.e., the character of $\bar\sigma$ is
  $-R_\finiteC^{\finiteG_y}\,\varphi$, where
  $\varphi=\varphi_1\otimes\varphi_2\otimes\varphi_3$ is a character
  of $\finiteC$ in general position.  Then, by~\cite[6.9]{srinivasan},
$$
\theta_{\bar\sigma} (\bar{\gamma}) =
    -\sum_{w\in W_{k_F}
    (\underline{\finiteC},\underline{\finiteG}_y)}
        {}^w(\varphi_1\otimes\varphi_2\otimes\varphi_3)
    (\bar{\gamma}).
$$
  An easy application of the character theory of abelian groups shows
  that there is some $\gamma$ of above type for which this sum does
  not vanish.  Thus if $\pi$ is of type \depthzerosc{2} or \depthzerosc{3},
  then $\pi$ is not stable.
Similarly, if $\pi$ is compactly induced from $G_z$ (i.e., $\pi$ is of
type \depthzerosc{4}),
then one can find stably conjugate $\gamma$ and $\gamma'$ such that
$\theta_\pi (\gamma)=0$, but $\theta_\pi
(\gamma')\neq 0$.  Hence $\pi$ is again not stable.
\end{proof}

\subsection{Non-singleton $L$-packets containing supercuspidals}
\label{sub:non-stable}

Let $\gamma\in C$ be a regular element of $G$ whose image $\bar\gamma$
in $\finiteC$ is a regular element of $\finiteG_y$.  Let $\gamma'\in
C'$ be the conjugate of $\gamma$ by $\nu$.  Since $\gamma\in G$ and
$\bar{\gamma}\in\finiteG$ are regular elliptic,
Lemma~\ref{lem:very_regular} implies that $\gamma$ lies in a
unique maximal parahoric subgroup, namely $G_y$.
Similarly, $\gamma'$ is not
contained in any parahoric subgroup other than $G_z$.
The following lemma then
follows easily from Proposition~\ref{induced_char},~\cite[6.9]{srinivasan},
and~\cite[p.~31]{ennola}.
\begin{lem}
\label{lem:char_value}
Let $\pi$ be a supercuspidal representation of $G$ of depth zero.
Then, in the notation of \S\ref{sec:depth-zero-G},
\begin{equation*}
\begin{aligned}
\theta_{\pi} (\gamma) &=
\begin{cases}
\displaystyle{
        -\sum_{w\in W_{k_F}(\underline{\finiteC},\underline{\finiteG}_y)}
 {}^w(\varphi_1\otimes\varphi_2\otimes\varphi_3) (\bar{\gamma})} 
 & \mbox{if $\pi$ is of type \textup{\depthzerosc{2}},}\\
2(\eta\otimes\eta\otimes\eta) (\bar{\gamma}) 
 & \mbox{if $\pi$ is of type \textup{\depthzerosc{3}},}\\
0 & \mbox{if $\pi$ is of type \textup{\depthzerosc{4}},}
\end{cases} \\
\theta_{\pi} (\gamma ') &= 
\begin{cases}
0 & \mbox{if $\pi$ is of type \textup{\depthzerosc{2}},}\\
0 & \mbox{if $\pi$ is of type \textup{\depthzerosc{3}},}\\
\displaystyle{-\sum_{w\in W_{k_F}(\underline{\finiteC}',\underline{\finiteG}_z)}
 {}^w(\varphi_1\otimes\varphi_2\otimes\varphi_3)
    (\bar{\gamma})}  
& \mbox{if $\pi$ is of type \textup{\depthzerosc{4}}.}
\end{cases}
\end{aligned}
\end{equation*}
\end{lem}

There are two types of non-singleton $L$-packets containing a
supercuspidal representation of $G$ of depth zero as discussed in
\S\ref{sec:depth-zero-G}; namely, the non-supercuspidal
$L$-packets of size two and the supercuspidal $L$-packets of
size four.  An $L$-packet $\Pi$ of the former type consists of
the unique square-integrable constituent $\pi^2 = \pi^2(\lambda)$
(see \S\ref{sec:depth-zero-G})
of a
reducible principal series of type~\reducible{2}
together with a corresponding
supercuspidal representation $\pi^s = \pi^s(\lambda)$.  Here $\lambda$
is a depth-zero character of $M$ such that $\lambda_1|_{F^\times} =
\omega_{E/F}|\cdot |_F^{\pm 1}$.  Recall the characters
$\lambda_1\in\Hom (E^\times ,\C^\times)$ and $\lambda_2\in\Hom (E^1
,\C^\times)$ determined by $\lambda$ according to~\eqref{eqn:lambda}.
Let $\bar\lambda_1$ and $\bar\lambda_2$ denote the associated
characters of $k_E^\times$ and $k_E^1$, respectively.

\begin{prop}
\label{prop:card2}
Let $\lambda$ be a depth-zero character of $M$ such that
$\lambda_1|_{F^\times} = \omega_{E/F}|\cdot |_F^{\pm 1}$.  Let
$\bar\lambda_1'$ denote the character of $k_E^1$ such that
$\bar\lambda_1'\circ\N = \bar\lambda_1$.
\begin{enumerate}
\item [(i)] If $\lambda_1$ is trivial on
$\OK_E^\times$, then 
$\pi^s (\lambda) = \ind_{G_y}^G \sigma,$
where $\sigma$ is the
inflation to $G_y$ of the representation 
$\tau\cdot (\bar\lambda_2\circ\det)$ 
of $\finiteG_y$.
\item [(ii)] If $\lambda_1$ is nontrivial on
$\OK_E^\times$, then 
$\pi^s (\lambda) = \ind_{G_z}^G \sigma,$
where $\sigma$ is the
inflation to $G_z$ of the representation of $\finiteG_z$ whose
character is 
$$
-R_{\finiteC'}^{\finiteG_z} (\bar\lambda_2 \otimes
\bar\lambda_1'\bar\lambda_2 \otimes \bar\lambda_1'\bar\lambda_2).
$$
\end{enumerate}
\end{prop}

\begin{proof}
  We determine the supercuspidal representation $\pi^s =
  \pi^s(\lambda)$ by computing its character at certain regular
  elliptic elements of $G$.  Recall that the irreducible constituents
  of $\ind_B^G\lambda$ are $\pi^2 (\lambda)$ and a non-tempered
  representation $\pi^n(\lambda)$ and that the set $\Pi' =
  \{\pi^s(\lambda),\pi^n(\lambda)\}$ is an $A$-packet of $G$.  Then
  $\Pi'$ is the endoscopic lift from $H$ to $G$ of the character
\begin{equation}
\label{eq:xi}
\xi = (\mu\lambda_2\circ\det_{\UU(1,1)}) \otimes \lambda_2
\end{equation}
of $H$, where $\mu\circ\N = \lambda_1 |\cdot |_E^{\mp
  1/2}\omega_{E'/E}$~\cite[\S12.2, \S13.1]{rogawski},
and $E'$ is an unramified quadratic extension of $E$.
Let $\omega$ be the central character of the elements of $\Pi$
and let $f\in C(G,\omega)$.  By~\cite[Thm.~13.1.1,
Prop.~13.1.2]{rogawski},
$$
\theta_{\pi^n}(f)+\theta_{\pi^s}(f) = \theta_\xi (f^H),
$$
where $f\mapsto f^H$ is the endoscopic transfer from $G$ to $H$
(see~\cite[\S4.3]{rogawski}).
Thus 
\begin{equation}
\label{eq:sc_char}
\theta_{\pi^s} = \theta_\xi^G - \theta_{\pi^n},
\end{equation}
where $\theta^G_\xi$ is the distribution on $G$ that arises from
$\theta_\xi$ via endoscopy.  The same equation holds for the functions
on $G\reg$ that represent these distributions.  Let $\gamma$ be an
element of $C$ whose image $\bar\gamma$ in $\finiteC$ is regular.  Let
$\gamma'\in C'$ be the conjugate of $\gamma$ by $\nu$ and let
$\bar\gamma'$ be its image in $\finiteC'$.  In order to determine
$\pi^s$, we will evaluate the right-hand side of (\ref{eq:sc_char}) at
$\gamma$ if $\lambda_1|_{\OK_E^\times}$ is trivial and at $\gamma'$ if
$\lambda_1|_{\OK_E^\times}$ is nontrivial.

First we compute $\theta_\xi^G(\gamma)$ and $\theta_\xi^G(\gamma')$.
By~\cite[Lemma 12.5.1]{rogawski} and the particular form of $\gamma$,
\begin{equation}
\label{eq:dist_transfer}
\theta_\xi^G(\gamma)
= \sum_{w\in W_F(\underline{C},\underline{H})
        \backslash W_F(\underline{C},\underline{G})}
                     \kappa(c_w)  \xi (^w\gamma),
\end{equation}
where $c_w$ is the class in
$$
\D (C/F) := \text{Ker}\{ H^1(F,\underline{C})\rightarrow
H^1(F,\underline{G})\}
$$
represented by the cocycle $\{ s(w)w^{-1}\}$
($s\in\Gal (\overline{F}/F)$) and $\kappa$ is the element of the dual
of $\D (C/F)$ corresponding to the endoscopic group $H$.  Since
$W_F(\underline{C},\underline{G})=W(C,G)$,
$\kappa (c_w)=1$ for all
$w\in W_F(\underline{C},\underline{H})\backslash
W_F(\underline{C},\underline{G})$.
Since $W_F(\underline{C},\underline{G})\cong S_3$
while $|W_F(\underline{C},\underline{H})|=2$, we
obtain
\begin{equation*}
\theta_\xi^G(\gamma) = \xi ((\gamma_1 ,\gamma_2,\gamma_3)) + \xi
((\gamma_3 ,\gamma_1,\gamma_2)) + \xi ((\gamma_2 ,\gamma_3,\gamma_1)).
\end{equation*}
Evaluating this 
when $\lambda_1|_{\OK_E^\times}$ is trivial and using (\ref{eq:xi}), we get
\begin{equation}
\label{eq:theta_gamma}
3 (\lambda_2\otimes\lambda_2\otimes\lambda_2)
(\gamma) = 3(\bar\lambda_2\otimes\bar\lambda_2\otimes\bar\lambda_2)
(\bar\gamma).
\end{equation}

As in the preceding paragraph,
$$
\theta_\xi^G(\gamma') = \sum_{w\in W_F(\underline{C}',\underline{H})
        \backslash W_F(\underline{C}',\underline{G})}
\kappa (c_w) \xi (^w\gamma') ,
$$
where $c_w$ is now the class in $\D (C'/F)$ represented by $\{ s(w)w^{-1}\}$.
In this case $W_F(\underline{C}',\underline{G})\cong S_3$
and $|W_F(\underline{C}',\underline{H})|=2$.  Then $\kappa (c_1) = 1$, and
an easy calculation shows that if $w$ represents a
nontrivial coset in
$W_F(\underline{C}',\underline{H})\backslash
W_F(\underline{C}',\underline{G})$, then $\kappa (c_w) = -1$ .
Thus 
\begin{equation*}
\theta_\xi^G (\gamma') = \xi((\gamma_1 ,\gamma_2,\gamma_3))
                     - \xi ((\gamma_3 ,\gamma_1,\gamma_2)) - \xi
                     ((\gamma_2 ,\gamma_3,\gamma_1)) .
\end{equation*}
We evaluate
this when $\lambda_1|_{\OK_E^\times}$ is
nontrivial.  Using (\ref{eq:xi}), we obtain
\begin{equation}
\label{eq:theta_gamma'}
(\bar\lambda_1'\bar\lambda_2 \otimes \bar\lambda_2 \otimes
\bar\lambda_1'\bar\lambda_2)
(\bar\gamma) -
(\bar\lambda_2 \otimes \bar\lambda_1'\bar\lambda_2 \otimes
\bar\lambda_1'\bar\lambda_2) 
(\bar\gamma) -
(\bar\lambda_1'\bar\lambda_2 \otimes \bar\lambda_1'\bar\lambda_2
\otimes \bar\lambda_2) 
(\bar\gamma).
\end{equation}

It remains to evaluate $\theta_{\pi^n}$ at $\gamma$ and $\gamma'$.
Since $\gamma$ and $\bar\gamma$ are regular elliptic and $\gamma\in
G_y$, $y$ is the unique fixed point of $\gamma$ in $\B$ by
Lemma~\ref{lem:very_regular}.  Then
\cite[Lemma~III.4.10, Theorem III.4.16]{schneider-stuhler} implies that 
\begin{equation}
\label{eq:schneider-stuhler}
\theta_{\pi^n}(\gamma) =
\text{trace}\left(\gamma|(\pi^n)^{G_{y+}}\right).
\end{equation}
The analogous formula holds for $\gamma'$ and $z$.  Hence we must
determine $(\pi^n)^{G_{y+}}$ and $(\pi^n)^{G_{z+}}$. 

Recall that $\pi^n(\lambda)$ and $\pi^2(\lambda)$ are the irreducible
constituents of $\ind_B^G\lambda$.  Let $\bar\lambda$ be the character
of $\finiteM$ determined by $\lambda$.  Since $G = G_yB$, we have that for
any $x\in\F$,
$$
\text{Res}_{G_y}\ind_{B}^{G}\lambda
=
\ind_{B\cap G_y} ^{G_y}\lambda
=
\ind_{G_x}^{G_y}
\ind_{B\cap G_y} ^{G_x}\lambda,
$$
which contains $\rho_y := \ind_{G_x} ^{G_y}\lambda $, the inflation to
$G_y$ of the representation $\bar\rho_y := \ind_{\finiteB_y}^{\finiteG_y}
\bar\lambda$.  Since $\rho_y$ is trivial on $G_{y+}$, this implies that
the space of $G_{y+}$-fixed vectors in $\ind_B^G\lambda$ contains
$\rho_y$.  Moreover, by Mackey's
theorem and Frobenius reciprocity,
\begin{align*}
\Hom_{G_{y+}} ({\bf 1}, \text{Res}_{G_{y+}}\ind_B^G\lambda )
&=
\Hom_{G_{y+}}
        \biggl({\bf 1}, 
        \bigoplus_{g\in G_{y+}\backslash G/B}
                \ind_{^gB\cap G_{y+}}^{G_{y+}}{}^g\!\lambda
        \biggr)\\ &= 
\bigoplus_{g\in G_{y+}\backslash G/B}
\Hom_{^gB\cap G_{y+}} \left({\bf 1}, {}^g\!\lambda \right)\\ &=
\bigoplus_{g\in G_{y+}\backslash G/B}
\Hom_{^gB\cap G_y+} \left({\bf 1}, {\bf 1}\right) .
\end{align*}
The dimension of this space is $|G_{y+}\backslash G/B|$, which (since
$G=G_yB$) is equal to 
$$
|G_{y+}\backslash G_y / (B\cap G_y)| = |\finiteG_y/\finiteB_y| =
\dim \rho_y .
$$
Hence the space of $G_{y+}$-fixed vectors in $\ind_B^G\lambda$ is
isomorphic to $\rho_y$.

Since the vertex $z$ is special, the Iwasawa decomposition
$G=G_z\overline{B}$ holds, where $\overline{B}$ is the Borel subgroup
opposite $B$ with respect to $M$.  Then an argument similar to that in
the preceding paragraph shows that, as a representation of
$\finiteG_z$, the space of $G_{z+}$-fixed vectors in
$\ind_{\overline{B}}^G\lambda$ is isomorphic to
$\bar\rho_z=\ind_{\finiteB_z}^{\finiteG_z}\bar\lambda$.

Now let $v$ equal $y$ if $\lambda_1|_{\OK_E^\times}$ is trivial or $z$
if $\lambda_1|_{\OK_E^\times}$ is nontrivial.  Let $\pi$ be either
$\pi^2$ or $\pi^n$.  By~\cite[Thm.~5.2]{moy-prasad2}, for $x\in\F$,
$(G_x , \lambda|_{M_0})$ is a $K$-type contained in $\pi$ (where we
have identified $G_x/G_{x+}$ and $M_0/M_{0+}$).  Thus, as a
representation of $\finiteB_v$, $\pi^{G_{x+}}$ contains the character
$\bar\lambda$ of $\finiteB_v$.  By Frobenius reciprocity,
$\pi^{G_{v+}}$ contains a subrepresentation of $\bar\rho_v$.  Since
$\bar\lambda$ extends to a character of $\finiteG_v$, $\bar\rho_v$ is
reducible with two irreducible constituents.  Replacing $\lambda$ by a
Weyl conjugate if necessary, we may assume that $\pi^2$ is a
subrepresentation of $\ind_B^G\lambda$, so that we have the exact
sequence
$$
0\longrightarrow \pi^2 \longrightarrow \ind_B^G\lambda
\longrightarrow \pi^n\longrightarrow 0.
$$
Taking $G_{v+}$-fixed vectors, we obtain the exact sequence
$$
0\longrightarrow (\pi^2)^{G_{v+}} \longrightarrow \bar\rho_v
\longrightarrow (\pi^n)^{G_{v+}}\longrightarrow 0
$$
of representations of $\finiteG_v$.  It follows that
as a representation of $\finiteG_v$,
$\pi^{G_{v+}}$ is an irreducible constituent of $\bar\rho_v$.

According to~\S\ref{sec:reps-finite}, the irreducible constituents of
$\bar\rho_v$ are a one-dimensional representation $\psi$ and the
representation $\St_{\finiteG_v} (\psi)$.  Here 
$$
\psi = (\bar\lambda'_1\circ\det_{\UU(1,1)}\circ p_v)\cdot
(\bar\lambda_2\circ\det_{\underline\finiteG_v}),
$$
where $p_v:
\underline\finiteG_v\longrightarrow U(1,1)$ is trivial if $v=y$ or the
projection onto the $U(1,1)$ factor of $\underline\finiteG_z\cong
U(1,1)\times U(1)$ if $v=z$.  Suppose that
$(\pi^2)^{G_{v+}}\cong\psi$.  Then $\finiteG_v$ acts via the character
$\psi$ on any nonzero vector $u\in (\pi^2)^{G_{v+}}$.  Let
$(\pi^2)\spcheck$ be the contragrediant representation of $\pi^2$.
Then $\finiteG_v$ acts via $\psi^{-1}$ on any nonzero vector
$u'\in((\pi^2)\spcheck)^{G_{v+}}$.  An easy computation shows that the
matrix coefficient $c_{u,u'}$ is not square-integrable.  It follows
that $u\notin\pi^2$ and hence that $(\pi^n)^{G_{v+}}\cong\psi$.  Thus,
if $\lambda_1|_{\OK_E^\times}$ is trivial, then from
\eqref{eq:schneider-stuhler},
$$
\theta_{\pi^n}(\gamma) =
\text{trace}\left(\gamma|(\pi^n)^{G_{y+}}\right) = \psi (\bar\gamma) =
(\bar\lambda_2\otimes\bar\lambda_2\otimes\bar\lambda_2 )(\bar\gamma).
$$
On the other hand, if $\lambda_1|_{\OK_E^\times}$ is nontrivial, then
from \eqref{eq:schneider-stuhler},
$$
\theta_{\pi^n}(\gamma ') =
\text{trace}\left(\gamma'|(\pi^n)^{G_{z+}}\right) = \psi (\bar\gamma') =
(\bar\lambda'_1\bar\lambda_2 \otimes \bar\lambda_2 \otimes
\bar\lambda'_1\bar\lambda_2 
)(\bar\gamma).
$$
Combining these calculations with (\ref{eq:theta_gamma}),
(\ref{eq:theta_gamma'}) and (\ref{eq:sc_char}), we find that if 
$\lambda_1|_{\OK_E^\times}$ is trivial,
\begin{equation}
\label{eq:theta_pi_s}
\theta_{\pi^s}(\gamma) =
2(\bar\lambda_2\otimes\bar\lambda_2\otimes\bar\lambda_2
)(\bar\gamma),
\end{equation}
while if $\lambda_1|_{\OK_E^\times}$ is nontrivial
\begin{equation*}
\theta_{\pi^s}(\gamma) = - (\bar\lambda_2 \otimes
\bar\lambda_1'\bar\lambda_2 \otimes
\bar\lambda_1'\bar\lambda_2) 
(\bar\gamma) -
(\bar\lambda_1'\bar\lambda_2 \otimes \bar\lambda_1'\bar\lambda_2
\otimes \bar\lambda_2) 
(\bar\gamma).
\end{equation*}

Suppose that $\lambda_1|_{\OK_E^\times}$ is trivial.  Since $\pi^s$ is
a depth-zero supercuspidal representation, Lemma~\ref{lem:char_value}
implies that $\theta_{\pi^s}(\gamma)$ is equal to the evaluation at
$\bar\gamma$ of a linear combination $\mu$ of characters of $\finiteC$
depending only on $\pi^s$.  Letting $\gamma$ vary over all elements of
$C$ that are regular in $G$ and that have regular image $\overline\gamma$
in $\finiteG_y$, we obtain from (\ref{eq:theta_pi_s}) that
$
\mu = 2(\bar\lambda_2\otimes\bar\lambda_2\otimes\bar\lambda_2)
$
on the set of regular elements of $\finiteC$.  By
Cor.~\ref{cor:q-value}, it must be the case that $\mu =
2(\bar\lambda_2\otimes\bar\lambda_2\otimes\bar\lambda_2)$.  By the
linear independence of characters of $\finiteC$, $\mu $ must have
the character $\bar\lambda_2\otimes\bar\lambda_2\otimes\bar\lambda_2$
as a summand.
Hence, by Lemma~\ref{lem:char_value}, $\pi^s$ must be equivalent to
$\ind_{G_y}^G \sigma,$ where $\sigma$ is the inflation to $G_y$ of
$\tau\cdot (\bar\lambda_2\circ\det)$.  This proves (i).  A similar
argument with $\gamma'$ replacing $\gamma$ proves (ii).
\end{proof}

We now determine the $L$-packets of $G$ of size 4.  Fix distinct
characters $\chi_1$, $\chi_2$, and $\chi_3$ of $k_E^1$.  Let $\chi$ be
the character $\chi_1\otimes\chi_2\otimes\chi_3$ of $k_E^1\times
k_E^1\times k_E^1$.  Define regular characters $\chi^{(1)}$,
$\chi^{(2)}$, and $\chi^{(3)}$ of $\finiteC'$ by
\begin{align*}
\chi^{(1)} &= \chi\\
\chi^{(2)} &= \chi_2\otimes\chi_3\otimes\chi_1\\
\chi^{(3)} &= \chi_3\otimes\chi_1\otimes\chi_2.
\end{align*}
Note that each $\chi^{(i)}$ is equal to $^w\chi$ for some $w\in
W_{k_F}(\underline{\finiteC}',\underline{\finiteG}_z)$.
Let $\sigma$ be the inflation to $G_y$
of the cuspidal representation $\bar\sigma$ of $\finiteG_y$ with
character $-R_\finiteC^{\finiteG_y}\,\chi$.  For $i=1,2,3$, let
$\sigma_i$ be the inflation to $G_z$ of the cuspidal representation
$\bar\sigma_i$ of $\finiteG_z$ with character
$-R_{\finiteC'}^{\finiteG_z}\,\chi^{(i)}$.  Then
$\sigma_1,\sigma_2,\sigma_3$ are distinct by \cite[p.~139]{srinivasan}.
Define $\pi_0 = \ind_{G_y}^G\sigma$ and
$\pi_i = \ind_{G_z}^G\sigma_i$ ($i=1,2,3$).
By~\cite{moy-prasad2}, these representations
are inequivalent supercuspidals of depth zero.  For $v=y$ or $z$, let
$\sigma_v$ be the inflation to $H_v$ of the cuspidal representation of
$\finiteH_v$ with character $-R_\finiteT^{\finiteH_v}\,\chi$, where
$\finiteT=\finiteC$ if $v=y$, and $\finiteT=\finiteC'$ if $v=z$.
Define $\rho_v = \ind_{H_v}^H \sigma_v$.  Then $\rho_y$ and $\rho_z$
are inequivalent but conjugate by an element of $\PGL_2(F)\times \{
1\}$, and hence $\{\rho_y,\rho_z\}$ is an $L$-packet for $H$.

\begin{prop}
\label{prop:card4}
The set $\{\pi_0,\pi_1,\pi_2,\pi_3\}$ is an $L$-packet for $G$ and is
the endoscopic transfer of $\{\rho ,\rho'\}$.
\end{prop}

\begin{proof}
  Let $R=\{\rho_y ,\rho_z\}$ and let $\Pi$ be the transfer of $R$ from
  $H$ to $G$.  Then $\Pi$ has size four
  by~\cite[Prop.~13.1.2]{rogawski}.  Let $\pi'_0,\pi'_1,\pi'_2,\pi'_3$
  be the elements of $\Pi$.  Then the $\pi'_i$ are supercuspidal
  by~\cite[Prop.~13.1.3(b)]{rogawski}.  That they have depth zero
  follows from our assumption (see the Introduction)
  that the transfer preserves depth.
  Set $\theta_R =
  \theta_{\rho_y} + \theta_{\rho_z}$.  Let $\theta_R^G$ be the
  endoscopic transfer of $\theta_R$ from $H$ to $G$.  It follows
  from~\cite[Thm.~13.1.1, Prop.~13.1.3, Lemma 12.7.2]{rogawski} that
  \begin{equation}
  \label{eq:card4}
  \theta_R^G = \theta_{\pi'_0} + \theta_{\pi'_1} - \theta_{\pi'_2} -
  \theta_{\pi'_3}
  \end{equation}
  for some ordering of the $\pi'_i$.
  Let $\gamma$ and $\gamma'$ be as in
  Proposition~\ref{prop:card2}.  We will compute $\theta_R(\gamma)$
  and $\theta_R(\gamma')$ to determine the $\pi'_i$.

  Let $\gamma^*$ be either $\gamma$ or $\gamma'$, and correspondingly
  let $T$ be either $C$ or $C'$.
  According to~\cite[Lemma 12.5.1]{rogawski}, using the notation in
  the proof of Proposition~\ref{prop:card2},
  $$
  \theta_R^G(\gamma^*) = \sum_{w\in W_F(\underline{T},\underline{H})\backslash
        W_F(\underline{T},\underline{G})}
  \kappa(c_w) \theta_R ({}^w\gamma^*).
  $$
  As in the proof of Proposition~\ref{prop:card2}, if $\gamma^* =
  \gamma$, then $\kappa (c_w)=1$ for all $w\in
  W_F(\underline{C},\underline{H})\backslash
  W_F(\underline{C},\underline{G})$, while if $\gamma^*=\gamma'$, then
  $\kappa (c_1) = 1$ and $\kappa (c_w) = -1$ if $w$ represents a
  nontrivial coset in $W_F(\underline{C}',\underline{H})\backslash
  W_F(\underline{C}',\underline{G})$.  Since $\gamma^*\in H$ and
  $\bar\gamma^*\in\finiteH$ are regular elliptic, $\gamma^*$ lies in a
  unique maximal parahoric subgroup $H_v$ of $H$ by
  Lemma~\ref{lem:very_regular} (where $v=y$ if $\gamma^* = \gamma$, and
  $v=z$ if $\gamma^* = \gamma'$).  Let $u$ be either $y$ or $z$.  It
  follows from Proposition~\ref{induced_char}
  and~\cite[6.9]{srinivasan} that
  $$
  \theta_{\rho_u} ({}^w\gamma^*) = 
  \left\{
    \begin{array}{ll}
      \displaystyle{
        -\sum_{u\in W_{k_F}(\underline{\finiteT},\underline{\finiteH})}
        {}^{uw}\chi
         (\bar{\gamma}^*)} &
         \mbox{if $u=v$}\\
       0 & \mbox{if $u\neq v$,}
    \end{array}
  \right.
  $$
  where we identify $W_{k_F}(\underline{\finiteT},\underline{\finiteH})$ with
$W_F(\underline{T},\underline{H})\subset W_F(\underline{T},\underline{G})$.
  Hence 
\begin{equation}
\label{eq:thetaR}
\begin{aligned}
  \theta^G_R (\gamma)
&= \displaystyle{-\sum_{w\in W_{k_F}(\underline{\finiteC},\underline{\finiteG})}
 {}^w\chi (\bar{\gamma}),}\\
  \theta^G_R (\gamma')
&= \displaystyle{-\sum_{w\in W_{k_F}(\underline{\finiteC},
                \underline{\finiteG})}d_w
  {}^w\chi (\bar{\gamma}),}
\end{aligned}
\end{equation}
where $d_w=1$ if $w\in W_{k_F}(\underline{\finiteC},\underline{\finiteH})$
and $d_w=-1$ otherwise.
 
 As observed in the proof of Proposition~\ref{prop:card4},
 Lemma~\ref{lem:char_value} implies that $\theta_{\pi'_i}(\gamma)$ is
 equal to the evaluation at $\bar\gamma$ of a linear combination
 $\mu_i$ of characters of $\finiteC$ depending only on $\pi'_i$.
 Therefore, evaluating (\ref{eq:card4}) at all $\gamma$ of the above
 type and using (\ref{eq:thetaR}), we obtain
$$
-\sum_{w\in W_{k_F}(\underline{\finiteC},\underline{\finiteG})}
        {}^w\chi = \mu_0+\mu_1-\mu_2-\mu_3
$$
 on the set of regular elements of
 $\finiteC$.  Then, by Cor.~\ref{cor:q-value}, this equation must hold
 at all elements of $\finiteC$.  It follows from
 Lemma~\ref{lem:char_value} and the linear independence of characters
 of $\finiteC$ that, after possibly reordering, $\pi'_0$ must be
 equivalent to $\pi_0$ and that the other elements of the $L$-packet
 must be induced from $G_z$.  
 Evaluating (\ref{eq:card4}) at
 $\gamma'$ of the above type and using a similar argument, we obtain that
 $\pi'_1\cong\pi_1$ and, up to reordering, $\pi'_i\cong\pi_i$ for
 $i=2,3$.
\end{proof}

\section{Explicit base change for $G$}
\label{sec:base_change}

\subsection{Packets consisting of principal series constituents}
\label{sub:bc_ps}

\begin{prop}
\label{prop:bc_ps}
Let $\lambda\in\Hom (M,\C^\times)$.
\begin{enumerate}[(i)]
\item
If $\ind_B^G\lambda$ is irreducible
  and $\ind_{\widetilde B}^{\widetilde G}\tilde\lambda$ is irreducible,
  then the
  base change lift of the $L$-packet $\{ \ind_B^G \lambda\}$ is
  $\ind_{\widetilde{B}}^{\widetilde{G}}\tilde\lambda$.
\item
If $\ind_B^G\lambda$ is irreducible
  but $\ind_{\widetilde B}^{\widetilde G}\tilde\lambda$ is reducible,
  then $\lambda_1|_{F^\times} = |\phantom{x}|_F^{\pm1}$, and
  the base change lift of
  $\{ \ind_B^G \lambda\}$ is
  $\ind_{\widetilde P}^{\widetilde G}
        \Bigl(
        ( \lambda_1\tilde\lambda_2|\cdot|_E^{\mp 1/2}\circ \det_{\GL(2)})
                \otimes \tilde\lambda_2
        \Bigr)$.
\item If $\lambda_1 = |\phantom{x} |_E^{\pm 1}$, then the lift of the
  $L$-packet comprising the one-dimensional constituent
  $\psi=\lambda_2\circ\det_{\underline{G}}$ (respectively, the
  Steinberg constituent $\St_G(\psi)$) of $\ind_B^G\lambda$ is the
  one-dimensional constituent $\tilde{\psi} =
  \tilde\lambda_2\circ\det_{\underline{\widetilde{G}}}$ (respectively,
  the Steinberg constituent $\St_{\widetilde{G}}(\tilde{\psi})$) of
  $\ind_{\widetilde{B}}^{\widetilde{G}}\tilde\lambda$.
\item  If $\lambda_1|_{F^\times}$ is trivial and $\lambda_1$ is
  nontrivial, then the lift of the $L$-packet
  $\{\pi_1(\lambda),\pi_2(\lambda)\}$ is
  $\ind_{\widetilde{B}}^{\widetilde{G}}\tilde\lambda$.
\end{enumerate}
\end{prop}

\begin{proof}
Cases~(i), (iii), and~(iv) 
follow from~\cite{rogawski} (Prop.~4.10.2 and the
paragraph before Theorem~13.2.1).
To prove case (ii), note that up
to the action of the Weyl group, we may assume that $\lambda$
is positive with respect to $B$.
The paragraph before Theorem~13.2.1 in~\cite{rogawski}
then implies that the base change lift of $\{ \ind_B^G \lambda\}$ is
the Langlands quotient of
$\ind_{\widetilde B}^{\widetilde G}\tilde\lambda$.
This quotient is the desired representation.
\end{proof}

\subsection{Stable supercuspidal representations}
\label{sub:bc_stable_sc}
\begin{figure}
$$
\begin{xy}
\xymatrix@+50pt{
\pi \ar@{|->}[r]^{\txt{base}}_{\txt{change}}
        \save+<-\labeldistance,\labeldistance>*{G} \restore
        &
\widetilde{\pi} \ar@{|->}[r]^{\txt{extend}}_{\txt{somehow}}
        \save+<0ex,\labeldistance>*{\widetilde{G}} \restore
        &
\widetilde{\pi}
        \save+<0ex,\labeldistance>*{\widetilde{G}\Gamma} \restore
        &
{\widetilde\pi'}
        \save+<\labeldistance,\labeldistance>*{\widetilde{G}\Gamma} \restore
\\
\sigma \ar@{|->}[u]^{\txt{Ind}}
        \save+<-\labeldistance,0cm>*{G_y} \restore
        &
&
\widetilde\sigma \ar@{|.>}[r]^{\txt{unique}}_{\txt{extension}}
        \save+<0ex,\labeldistance>*{\widetilde{G}_y\Gamma} \restore
        &
{\widetilde\sigma}
        \ar@{|->}[u]_{\txt{Ind}}
        \save+<\labeldistance,-\labeldistance>*
        {\widetilde Z \widetilde G_y \Gamma}\restore
\\
\overline\sigma \ar@{|->}[r]^{\txt{base}}_{\txt{change}}
        \ar@{|->}[u]^{\txt{inflate}}
        \save+<-\labeldistance,-\labeldistance>*{\finiteG} \restore
        &
\widetilde{\overline\sigma} \ar@{|->}[r]^{\txt{extend}}_{\txt{somehow}}
        \save+<0ex,-\labeldistance>*{\widetilde{\finiteG}} \restore
        &
\widetilde{\overline\sigma} \ar@{|.>}[u]^{\txt{inflate}}
                                \ar@{|->}[ur]_{\txt{inflate}}
        \save+<0ex,-\labeldistance>*{\widetilde{\finiteG}\Gamma} \restore
}
\end{xy}
$$
\caption{Two ways of constructing representations of $\widetilde{G}\Gamma$}
\label{fig:stable}
\end{figure}
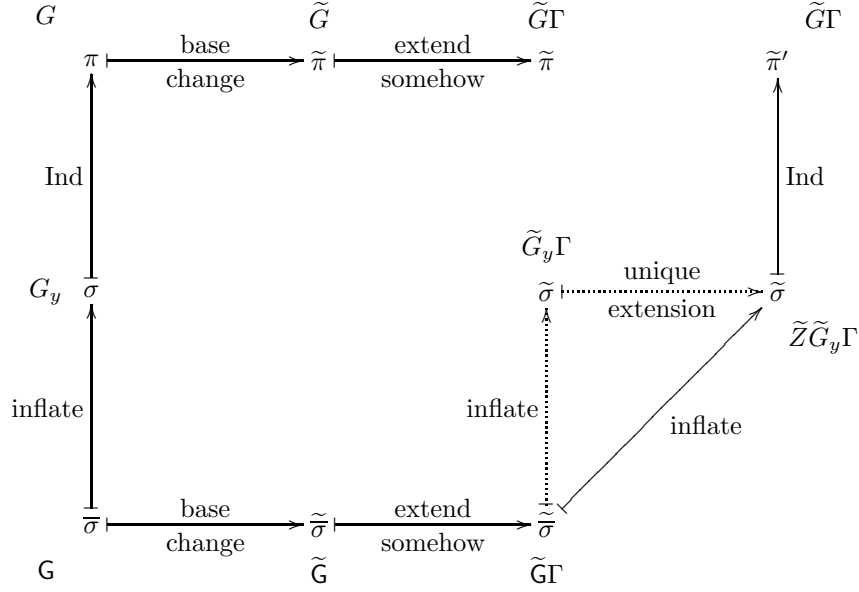

Suppose $\pi$ is a depth-zero, stable, supercuspidal representation of $G$.
From Prop.~\ref{sub:stable_sc}, $\pi^{G_y^+}$ contains
the inflation $\sigma$ of a cubic cuspidal
representation $\bar\sigma$ of $\finiteG\cong\finiteG_y$.
Then Figure~\ref{fig:stable} illustrates how to construct
representations $\tilde\pi$ and $\tilde\pi'$ of
$\widetilde{G}\Gamma$.
We can describe base change for $\pi$ explicitly
by showing that
$\tilde\pi$ and $\tilde\pi'$ are equivalent, provided
that the extensions from $\widetilde{G}$ to $\widetilde{G}\Gamma$
and from $\widetilde\finiteG$ to $\widetilde\finiteG\Gamma$
are chosen in compatible ways.

\begin{rem}
\label{rem:cartan}
Recall the Cartan decomposition for $\widetilde{G}$:
The diagonal subgroup $\widetilde{M}$ determines
a root system $\Phi$ for $\widetilde{G}$, and
the Borel subgroup $\widetilde{B}$ determines a positive
root system $\Phi^+$ inside $\Phi$.
Let $\widetilde{M}^+$ denote the set of all $m\in\widetilde{M}$ 
such that $\alpha(m)$ has positive valuation for all $\alpha\in\Phi^+$.
Then
$$
\widetilde{G} = \bigcup_{m\in\widetilde{M}^+}\widetilde{G}_y m
 \widetilde{G}_y .
$$
 Moreover, $m,m'\in \widetilde{M}^+$ represent the same double coset
 if and only if $m'\in m\widetilde{M}_0$.
\end{rem}
\begin{lem}
\label{lem:coset_rep}
Every conjugate of $\widetilde{Z}\widetilde{G}_y\epsgen$ in
$\widetilde{G}\epsgen$ is of the form
$^{gm}\bigl(\widetilde{Z}\widetilde{G}_y\epsgen\bigr)$, where
$g\in\widetilde{G}_y$, $m\in\widetilde{M}^+$.
\end{lem}
\begin{proof}
  The normalizer of $\widetilde{Z}\widetilde{G}_y\epsgen$ in
  $\widetilde{G}\epsgen$ is $\widetilde{Z}\widetilde{G}_y\epsgen$
  itself.  Therefore, the conjugates of $\widetilde{Z}\widetilde{G}_y\epsgen$
  correspond to the cosets in
  $\widetilde{G}\epsgen/\widetilde{Z}\widetilde{G}_y\epsgen\cong
  \widetilde{G}/\widetilde{Z}\widetilde{G}_y$.  The lemma now follows
  from the Cartan decomposition.
\end{proof}

\begin{rem}
\label{rem:vertex-types}
Recall that an inner automorphism of $\widetilde{G}$ acts
on the extended Dynkin diagram either trivially or via
a rotation.  Thus, if $g\in\widetilde{G}$ stabilizes the alcove
$\Delta\subset\widetilde{\B}^{\text{red}}$ and fixes some point
in the closure of $\Delta$ not equal to the barycenter of $\Delta$,
then $g$ must fix $\Delta$ pointwise.
\end{rem}

In the next result, we use the fact that $\widetilde{\finiteG}_y$ (in
addition to being a quotient of $\widetilde{G}_y$) is a quotient of
$\widetilde{Z}\widetilde{G}_y$.

\begin{prop}
\label{prop:bc_stable}
Let $\sigma$ be the inflation to $G_y$ of a cubic cuspidal
representation $\bar\sigma$ of $\finiteG_y$, and let $\tilde{\bar\sigma}$
be the Shintani lift of $\bar\sigma$ from $\finiteG_y$ to
$\widetilde{\finiteG}_y$.  Let $\pi=\ind_{G_y}^G\sigma$.  Then the
base change lift of the $L$-packet $\{\pi\}$ is 
$$
\ind_{\widetilde{Z}\widetilde{G}_y}^{\widetilde{G}}\tilde{\sigma},
$$
where $\tilde{\sigma}$ is the inflation to
$\widetilde{Z}\widetilde{G}_y$ of $\tilde{\bar\sigma}$.
\end{prop}

\begin{proof}
Let
$\tilde{\pi}$ be the base change lift of $\{\pi\}$
and let
$\ptp =
        \ind_{\widetilde{Z}\widetilde{G}_y}^{\widetilde{G}}\tilde{\sigma}$.
Since
$\theta_\pi$ is stable by Proposition~\ref{prop:stable_sc},
$\theta_{\tilde{\pi},\varepsilon}$
is a stable $\varepsilon$-class function on
$\widetilde{G}\ereg$ according to~\cite[\S12.5]{rogawski}.
By~\cite[Prop.~6.8]{moy-prasad2},
$\ptp$ is a supercuspidal representation
of $\widetilde{G}$ of depth zero.  Also, 
$$
{}^\varepsilon\ptp \cong
\ind_{\widetilde{Z}\widetilde{G}_y}^{\widetilde{G}}{}
^\varepsilon\tilde{\sigma} \cong \ptp
$$
since $\tilde{\sigma}$ is $\varepsilon$-invariant as it is in the
image of the Shintani lift.  Since $\ptp$ is $\varepsilon$-invariant,
$\ptp$ is the base change lift of a singleton supercuspidal $L$-packet
$\{\pi'\}$ by~\cite[Prop.~13.2.2]{rogawski}.  But then, as in the case
of $\tilde{\pi}$, $\theta_{\ptp,\varepsilon}$ is a stable
$\varepsilon$-class function.  Moreover, it is easily seen (under the
assumption that $\pi$ has depth zero) that the central characters of
$\tilde{\pi}$ and $\ptp$ are identical.  Furthermore, according to
~\cite[\S13.2]{rogawski}, we may choose $\tilde\pi (\varepsilon)$ and
$\ptp (\varepsilon)$ so that $\theta_{\tilde{\pi},\varepsilon} =
\theta_{\pi}\circ\N$ and $\theta_{\ptp ,\varepsilon} =
\theta_{\pi'}\circ\N$ (see \S\ref{sec:notation-general}).

If $\alpha_1$ and $\alpha_2$ are stable $\varepsilon$-class functions
on $\widetilde{G}\ereg$ that transform under $\widetilde{Z}$ via the
same character, then the $\varepsilon$-elliptic inner product of
$\alpha_1$ and $\alpha_2$ (see~\cite[\S12.5]{rogawski}) is defined by
\begin{equation}
\label{eq:ell_inn_prod}
\langle\alpha_1,\alpha_2\rangle_\varepsilon =
 \sum_{T\in\mathcal{C}}|W_F(\underline{T},\underline{G})|^{-1}
 \int_{\widetilde{Z}\widetilde{T}^\N\backslash\widetilde{T}}
 D_G(\N(\delta))^2\alpha_1(\delta)\overline{\alpha_2(\delta)}\,d\delta ,
\end{equation}
where $\mathcal{C}$ is a set of representatives for the stable
conjugacy classes of elliptic Cartan subgroups of $G$, $\widetilde{T}$
is the centralizer of the Cartan subgroup $T$ in $\widetilde{G}$,
$\widetilde{T}^\N$ is the kernel of the norm map on $\widetilde{T}$,
and $D_G$ is the discriminant.

By~\cite[Prop.~12.6.2]{rogawski}, in order to prove that
$\tilde{\pi}\cong\ptp $, it suffices to show that
$$
\langle
\theta_{\ptp ,\varepsilon},\theta_{\tilde{\pi},\varepsilon}
\rangle_\varepsilon \neq 0.
$$
We will verify the non-vanishing of this inner product by showing
that the two twisted characters agree on $\widetilde{T}$ for each
$T\in\mathcal{C}$.

By Proposition~\ref{prop:stable_sc}, the stability of $\pi'$ implies
that $\pi'$, like $\pi$, is induced from the inflation $\sigma'$ to
$G_y$ of a cubic cuspidal representation $\bar\sigma'$ of
$\finiteG_y$.  Since $\bar\sigma$ and $\bar\sigma'$ both arise via
Deligne-Lusztig induction from a cubic Cartan subgroup of
$\finiteG_y$, these representations agree on unipotent elements of
$\finiteG_y$~\cite[6.9]{srinivasan}.  It follows that $\theta_\pi$ and
$\theta_{\pi'}$ agree on $G_{0+}$.  Since $\pi$ and $\pi'$ have the
same central character, $\theta_\pi$ and $\theta_{\pi'}$ agree on
$ZG_{0+}$.

Let $T\in\mathcal{C}$,
$\delta\in\widetilde{T}\cap\widetilde{G}\ereg$,
and $\gamma = \N(\delta)$.
If $\gamma\in ZT_{0+}$, then
$$
\theta_{\tilde{\pi},\varepsilon}(\delta)
= \theta_{\pi}(\gamma )
= \theta_{\pi'}(\gamma )
= \theta_{\ptp,\varepsilon}(\delta).
$$
  We may therefore assume that $\gamma\in T\smallsetminus ZT_{0+}$.

If no conjugate of $\gamma$ is contained in $G_y$, then both
$\theta_{\tilde{\pi},\varepsilon}(\delta) = \theta_{\pi}(\gamma )$ and
$\theta_{\ptp,\varepsilon}(\delta) = \theta_{\pi'}(\gamma )$
vanish by Proposition~\ref{induced_char}.
We may therefore assume that $\gamma\in G_y$.

We may assume that that $T$ is not of type \cartan{0}, since
such tori are not elliptic.

Suppose that $T$ is of type \cartan{1} or \cartan{2}.
As in the proof of
Proposition~\ref{prop:stable_sc}, since
$\gamma\in T\smallsetminus ZT_{0+}$,
the semisimple part of the image $\bar{\gamma}$ of $\gamma$ in $\finiteG_y$
is not contained in a
cubic torus of $\finiteG_y$.  Therefore, by~\cite[6.9]{srinivasan}, 
$\theta_{\bar\sigma}$ and
$\theta_{\bar\sigma'}$ vanish on $\gamma$.  Thus, again we have
$$
\theta_{\tilde{\pi},\varepsilon}(\delta) = \theta_{\pi}(\gamma ) = 0 =
\theta_{\pi'}(\gamma )
=
\theta_{\ptp,\varepsilon}(\delta)
$$
by Proposition~\ref{induced_char}.

Now suppose that $T$ is of type \cartan{3}.  Then there exist cubic
extensions $L$ of $E$ and $K$ of $F$ such that $L=EK$ and
$T\cong\text{Ker}(N_{L/K})$.  In particular, we may identify $T_{0+}$
with $\text{Ker}(N_{L/K})\cap (1+\mathfrak{p}_L)$ and $Z$ with $E^1$.
We therefore have $ZT_{0+}\cong E^1\left[\text{Ker}(N_{L/K})\cap
  (1+\mathfrak{p}_L)\right]$.  If $L/E$ is totally ramified, then
$E^1\left[\text{Ker}(N_{L/K})\cap (1+\mathfrak{p}_L)\right] =
\text{Ker}(N_{L/K})$ so $T\smallsetminus ZT_{0+}$ is empty, and there
is nothing to prove in this case.  We may hence assume that $L/E$ is
unramified.  Since $T$ is determined only up to stable conjugacy, we
may also assume that $T$ fixes the point $y$.

Let $\bar{\gamma}$ be the image of $\gamma$ in
the cubic torus
$\finiteT\subset\finiteG_y$.  Since $\gamma\notin ZT_{0+}$,
$\bar{\gamma}$ is not central in $\finiteG_y$.  Thus $\bar{\gamma}$ is
regular elliptic, so by Lemma~\ref{lem:very_regular}, $\gamma$ is
contained in a unique parahoric subgroup of $G$, namely $G_y$.  Thus
$$
\theta_{\tilde{\pi},\varepsilon}(\delta) = \theta_{\pi}(\gamma) =
\theta_{\bar\sigma}(\bar{\gamma})
$$
by Proposition~\ref{induced_char}.  It
suffices to show that $\theta_{\ptp,\varepsilon}(\delta) =
\theta_{\bar\sigma}(\bar{\gamma}).$

Extend $\ptp$ to a representation (also denoted $\ptp$) of
$\widetilde{G}\epsgen$ in a manner compatible with the choice of
$\ptp(\varepsilon)$ made in the beginning of the proof.  Then
$$
\theta_{\ptp,\varepsilon}(\delta) =
\theta_{\ptp}(\delta\varepsilon).
$$
As a representation of $\widetilde{G}\epsgen$, 
$$
\ptp \cong
\ind_{\widetilde{Z}\widetilde{G}_y\epsgen}
    ^{\widetilde{G}\epsgen}
    \,\tilde\sigma,
$$
where $\tilde{\sigma}$ is extended compatibly from
$\widetilde{G}_y$ to $\widetilde{G}_y\epsgen$.  
This extension determines an extension of
$\tilde{\bar\sigma}$ to $\widetilde\finiteG_y\epsgen$, and we let
$\theta_{\bar\sigma}$ be the corresponding twisted character.

By Proposition~\ref{induced_char}, to compute
$\theta_{\ptp}(\delta\varepsilon)$, we must determine which
conjugates of $\widetilde{Z}\widetilde{G}_y\epsgen$
contain $\delta\varepsilon$.  Since
$$
\widetilde{T}^\N\widetilde{T}_0 \cong K^\times\OK_L^\times =
L^\times \cong \widetilde{T},
$$
and since
$\theta_{\ptp,\varepsilon}(\delta)$ only depends on $\delta$ modulo
$\widetilde{T}^\N$, we may assume that $\delta\in\widetilde{T}_0$.
Since $T$ fixes $y$, $\delta\in\widetilde{T}_0\subset \widetilde{G}_y$ so
$\delta\varepsilon\in\widetilde{G}_y\epsgen$.

Now suppose that $\delta\varepsilon$ is also contained in another
conjugate of $\widetilde{Z}\widetilde{G}_y\Gamma$.  By
Lemma~\ref{lem:coset_rep}, any such conjugate is of the form
$^{gm}(\widetilde{Z}\widetilde{G}_y\epsgen)$, where $g\in\widetilde{G}_y$ and
$m\neq 1$ is in $\widetilde{M}^+$.  If $\delta\varepsilon\in
{}^{gm}(\widetilde{Z}\widetilde{G}_y\epsgen)$, then
$^{(gm)^{-1}}(\delta\varepsilon)\in \widetilde{Z}\widetilde{G}_y\epsgen$ so
$$
m^{-1}\delta'\varepsilon (m)\in \widetilde{Z}\widetilde{G}_y,
$$
where $\delta'=g^{-1}\delta\varepsilon(g)\in\widetilde{G}_y$.  Thus
$\delta'\varepsilon (m)\in m\widetilde{Z}\widetilde{G}_y$ so
$\varepsilon (m)$ and $mc$ represent the same double coset in
$\widetilde{G}_y\backslash\widetilde{G}/\widetilde{G}_y$ for some
$c\in\widetilde{Z}$.  Since $\varepsilon(m),mc\in\widetilde{M}^+$, we
have $\varepsilon(m)\in mc\widetilde{M}_0$ by Remark~\ref{rem:cartan}.
It follows that $\delta'$ is in ${}^m(\widetilde{Z}\widetilde{G}_y) =
\widetilde{G}_{my}\widetilde{Z}$ as well as $\widetilde{G}_y$.  

Let $\bar{y}$ be the image of $y$ in $\widetilde{\B}^{\text{red}}$.
Then $\delta'$ fixes $\bar{y}$ and $m\bar{y}$, hence fixes the line
segment $[\bar{y},m\bar{y}]$ in $\widetilde\B^{\text{red}}$.  Since
$\varepsilon (m)\in mc\widetilde{M}_0$, we have $\varepsilon
(m\bar{y}) = m\bar{y}$ so that $[\bar{y},m\bar{y}]$ intersects the
(open) triangle $\Delta\subset\widetilde{\B}^{\text{red}}$ non-trivially.
Hence $\delta'$ stabilizes $\Delta$.  From
Remark~\ref{rem:vertex-types}, $\delta'$ must fix $\Delta$ pointwise.
Thus $\delta'$ is contained in
$\widetilde{G}_{\widetilde\F}\widetilde{Z}$, where
$\widetilde{G}_{\widetilde\F}$ is the standard upper-triangular
Iwahori subgroup of $\widetilde{G}$.  The image $\bar{\delta}'$ of
$\delta'$ in $\widetilde{\finiteG}_y$ is therefore contained in the
Borel subgroup $\widetilde{\finiteB}_y$ of upper-triangular matrices
in $\widetilde{\finiteG}_y$.  Since $\widetilde{\finiteB}_y$ is
$\varepsilon$-invariant, $\bar{\delta}'\varepsilon(\bar{\delta}')$ is
also contained in $\widetilde{\finiteB}_y$.  Hence the eigenvalues of
$\bar{\delta}'\varepsilon(\bar{\delta}')$ lie in $k_E^\times$.  But
$\bar{\delta}'\varepsilon(\bar{\delta}') =\bar\N(\bar\delta')
=\bar\N(\bar g^{-1} \bar\delta \varepsilon(\bar g)) =\bar g^{-1}
\bar\gamma \bar g$, where $\bar g$ is the image of $g$ in
$\widetilde\finiteG_y$, so the eigenvalues of
$\bar{\delta}'\varepsilon(\bar{\delta}')$ are the same as those of
$\bar{\gamma}$.  The eigenvalues of $\bar{\gamma}$, however, lie in
$k_L^1 \smallsetminus k_E^1$ since $\bar{\gamma}$ is a regular element
of the cubic torus $\finiteT$.  This contradiction shows that
$\delta\varepsilon$ is contained in a unique conjugate of
$\widetilde{Z}\widetilde{G}_y\epsgen$, namely
$\widetilde{Z}\widetilde{G}_y\epsgen$ itself.

We therefore have from Proposition~\ref{induced_char} that
$$
\theta_{\ptp,\varepsilon}(\delta) =
\theta_{\ptp}(\delta\varepsilon) =
\theta_{\tilde{\bar\sigma}}(\bar{\delta}\varepsilon)
= \theta_{\tilde{\bar\sigma},\varepsilon}(\bar{\delta}).
$$
But $\tilde{\bar\sigma}$ is the Shintani lift of $\bar\sigma$
(see~\cite{kawanaka}) so the last
expression is equal to 
$$
\pm\theta_{\bar\sigma}(\bar{\gamma}).
$$
(Here the twisted character
$\theta_{\tilde{\bar\sigma},\varepsilon}$ as chosen above is not a
priori equal to $\theta_{\bar\sigma}\circ\N$ since this choice is not
necessarily the one that is compatible with the Shintani lifting.
Nevertheless, it is at worst off by a sign by the discussion
in~\ref{sec:notation-general}.)  At the same time
$$
\theta_{\ptp,\varepsilon}(\delta) = \theta_{\pi'}(\gamma) =
\theta_{\bar\sigma'}(\bar{\gamma}),
$$
so $\theta_{\bar\sigma'}(\bar{\gamma}) =
\pm\theta_{\bar\sigma}(\bar{\gamma})$.  It is easily seen (e.g., from
the character table in~\cite{ennola}) that there is no cubic
cuspidal representation $\bar\sigma'$ of $\finiteG_y$ satisfying
$\theta_{\bar\sigma'}(\bar{\gamma}) =
-\theta_{\bar\sigma}(\bar{\gamma})$ for all regular elements
$\bar{\gamma}$ of cubic tori.  Thus
$$
\theta_{\ptp,\varepsilon}(\delta) =
    \theta_{\bar\sigma}(\bar{\gamma}),
$$
and the theorem follows.
\end{proof}

\subsection{Non-singleton $L$-packets containing supercuspidals}
\label{sub:bc_non-stable}

\begin{prop}
\label{prop:bc_card_2}
Let $\lambda$ be a character of $M$ of depth zero
such that $\lambda_1|_{F^\times}
= \omega_{E/F}|\cdot |_F^{\pm 1}$.
\begin{enumerate}[(i)]
\item The base change lift of the $L$-packet
$\{ \pi^2 (\lambda),\pi^s (\lambda)\}$ is 
$$
\ind_{\widetilde{P}}^{\widetilde{G}} \left(\St_{\widetilde{H}} \left( 
(\lambda_1 \tilde{\lambda}_2 |\cdot |_E^{\mp 1/2} \circ\det_{\GL(2)})
\otimes \tilde{\lambda}_2  
\right)\right).
$$
\item The base change lift of the $A$-packet
$\{ \pi^n (\lambda),\pi^s (\lambda)\}$ is 
$$
\ind_{\widetilde{P}}^{\widetilde{G}} \left( 
(\lambda_1 \tilde{\lambda}_2 |\cdot |_E^{\mp 1/2} \circ\det_{\GL(2)} )
\otimes \tilde{\lambda}_2  
\right).
$$
\end{enumerate}
Moreover, the above two base change lifts are precisely the
irreducible constituents of the principal series representation
$\ind_{\widetilde{B}}^{\widetilde{G}} (\tilde\lambda)$.
\end{prop}

Note that the proposition has the same content if we restrict
the choice of exponent in the hypothesis to be $+1$
(or to be $-1$).

\begin{proof}
This follows from~\cite[\S\S12--13]{rogawski}.  More
precisely, let $\xi$ be the character 
$$
(\mu\lambda_2\circ\det_{\UU(1,1)})
\otimes \lambda_2
$$
 of $H$, where 
$$
\mu\circ\N = \lambda_1 |\cdot
|_E^{\mp 1/2}\omega_{E'/E},
$$
$E'$ an unramified quadratic extension of $E$.
(Here, we are identifying $H$ with $\UU(1,1)(F)\times \UU(1)(F)$.)
  Let $\rho = \St_H (\xi)$.  Then, by
\cite[Prop.~13.1.3(c)]{rogawski}, the $L$-packet
$\{ \pi^2 (\lambda),\pi^s (\lambda)\}$ on $G$ is the lift of the $L$-packet
$\{\rho\}$ on $H$.  It follows from \cite[Prop.~13.2.2
(c)]{rogawski} that the base change lift of
$\{ \pi^2 (\lambda),\pi^s (\lambda)\}$
is $\ind_{\widetilde{P}}^{\widetilde{G}}(\tilde{\rho}')$,
where $\tilde{\rho}'$ is the ``primed'' base change lift
(see~\cite[\S11.4]{rogawski})
of $\rho$
from $H$ to $\widetilde{H}$.  But by \cite[\S12.1]{rogawski}, 
$$
\tilde{\rho}' = \St_{\widetilde{H}}(\tilde{\xi}'),
$$
where $\tilde{\xi}'$ is the character 
$$
(\omega_{E'/E}(\mu\lambda_2\circ\N )\circ\det_{\GL(2)})\otimes
\tilde\lambda_2 = (\lambda_1 \tilde{\lambda}_2 |\cdot |_E^{\mp 1/2}
\circ\det_{\GL(2)} ) \otimes \tilde{\lambda}_2.
$$
This proves (i), and
(ii) follows analogously from~\cite[Prop.~13.1.3(d)]{rogawski}.

The final statement follows from the
proof of~\cite[Lemma~12.7.6]{rogawski}.
\end{proof}

Recall the notation of Proposition~\ref{prop:card4}.  Let $\Pi$ be the
supercuspidal $L$-packet $\{\pi_0,\pi_1,\pi_2,\pi_3\}$,
and let $R=\{\rho_y,\rho_z\}$ be the $L$-packet of $H$ that
transfers to $\Pi$.
\begin{prop}
\label{prop:bc_card4}
The base change lift of the $L$-packet $\Pi = \{ \pi_0
,\pi_1,\pi_2,\pi_3\}$ is
$\ind_{\widetilde{B}}^{\widetilde{G}}\chi^*$, where $\chi^*$ is
inflation to $\widetilde{M}\cong E^\times\times E^\times\times
E^\times$ of
$$
\hat{\chi} = \tilde\chi_1\otimes
\tilde\chi_2 \otimes \tilde\chi_3
\in\Hom(\widetilde\finiteM,\C^\times).
$$
\end{prop}

\begin{proof}
  This also follows from~\cite[\S\S12--13]{rogawski}.  Note that since
  $E/F$ is unramified, $\widetilde{\finiteM}$ is a quotient of
  $\widetilde{M}$, so the definition of $\chi^*$ makes sense.  Let
  $\tilde{\rho}'$ be the ``primed'' base change lift
(see~\cite[\S11.4]{rogawski})
of $R$ from $H$
  to $\widetilde{H}$.  By \cite[Prop.~13.2.2(c)]{rogawski}, the base
  change lift of $\Pi$ is
  $\ind_{\widetilde{P}}^{\widetilde{G}}(\tilde{\rho}')$.
  By~\cite[\S12.1]{rogawski}, $R$ is the transfer from $C$ to $H$ of
  some character $\varphi $ of $C$.  Let $\theta^H_\varphi$ be the
  distribution on $H$ that arises from $\theta_\varphi = \varphi$
via endoscopy.
  Let $\widetilde{B}'$ be a Borel subgroup of $\widetilde{H}$
  containing $\widetilde{C}$.  Then by \cite[\S12.1]{rogawski},
  $\tilde{\rho}'$ is the representation
  $\ind_{\widetilde{B}'}^{\widetilde{H}}\tilde{\varphi}$.  Hence the
  base change lift of $\Pi$ is
  $$
  \ind_{\widetilde{P}}^{\widetilde{G}}
  \ind_{\widetilde{B}'}^{\widetilde{H}}\tilde{\varphi} =
  \ind_{\widetilde{B}'}^{\widetilde{G}}\tilde{\varphi}.
  $$
  We now determine $\varphi$.

  Since $R$ has depth zero, $\varphi$ must as well.  Also,
  \begin{equation}
  \label{eq:theta_phi}
  \theta_\varphi^H = \pm (\theta_\rho - \theta_{\rho'})
  \end{equation}
  by~\cite[Prop.~11.1.1(b)]{rogawski}.  The same equation holds for
  the functions that represent these distributions.
  Let $\gamma$ be an element of $C$ whose image $\bar\gamma\in \finiteC$
  is regular in $\finiteH_y$.  As computed in the proof of
  Propositon~\ref{prop:card4},
  $$
  \theta_\rho (\gamma) =
        -\sum_{w\in W_{k_F}(\underline{\finiteC} ,\underline{\finiteH})}
        {}^w\chi(\bar\gamma),
  $$
  while $\theta_{\rho'} (\gamma) = 0$.
  Hence the evaluation of the right side of (\ref{eq:theta_phi}) at
  $\gamma$ is 
  $$
  \pm\sum_{w\in W_{k_F}(\underline{\finiteC} ,\underline{\finiteH})}
        {}^w\chi(\bar\gamma).
  $$
  The analogue of (\ref{eq:dist_transfer}) for the transfer from
  $C$ to $H$ implies that
  $$
  \theta_\varphi^H (\gamma) =
        \sum_{w\in W_F(\underline{C},\underline{H})}\varphi ({}^w\gamma) =
  \sum_{w\in W_{k_F}(\underline{\finiteC} ,\underline{\finiteH})}
        {}^w\bar\varphi(\bar\gamma),
  $$
  where $\bar\varphi$ is the character of $\finiteC$ determined by
  $\varphi$.  Using (\ref{eq:theta_phi})
and letting $\gamma$ vary
  over all elements of the above type, it follows that
  $$
  \sum_{w\in W_{k_F}(\underline{\finiteC} ,\underline{\finiteH})}
        {}^w\bar\varphi =
  \pm\sum_{w\in W_{k_F}(\underline{\finiteC} ,\underline{\finiteH})}{}^w\chi
  $$
  on the set of regular elements of $\finiteC$, hence on all of
  $\finiteC$ by Cor.~\ref{cor:q-value}.  By the linear independence of
  characters of $\finiteC$, it follows that $\bar\varphi = {}^w\chi$
  for some $w\in W_{k_F}(\underline{\finiteC} ,\underline{\finiteG})$.
  Since $\tilde\varphi$
  is in the image of the base change lifting from $C$ to $\widetilde{C}$,
  it follows from~\cite[\S12.4]{rogawski}
  that $\tilde\varphi$ is trivial on elements of $\widetilde{C}$ of the
  form $(\varpi^a,\varpi^b,\varpi^c)$.  Since $\tilde\varphi$ has
  depth zero, $\tilde\varphi$ must be
  the inflation to $\widetilde{C}$ of $^w\tilde{\chi}$ for some $w$.
  Thus
  $^{w^{-1}}\tilde{\varphi}$ is the inflation to
  $\widetilde{C}$ of $\tilde\chi$, where $w$ is viewed as an element of
  $W(\widetilde{C},\widetilde{G})$.  Moreover,
$$
\tilde\pi = \ind_{\widetilde{B}'}^{\widetilde{G}}\tilde\varphi \cong
\ind_{\widetilde{B}'}^{\widetilde{G}}{} ^{w^{-1}}\tilde\varphi.
$$

Finally, note that by conjugating by a suitable element, one can send
$\widetilde{B}'$, $\widetilde{C}$, and $^{w^{-1}}\tilde{\varphi}$
respectively to $\widetilde{B}$, $\widetilde{M}$, and $\chi^*$.
The theorem follows.
\end{proof}

\section{Compatibility of base change and $K$-types}
\label{sec:K-types}

In this section we prove the Main Theorem, as stated
in \S\ref{sec:intro}.  Throughout, $\Pi$ will denote an $L$-packet of
$G$ and $\tilde\pi$ the base change lift of $\Pi$.

\subsection{Principal series $L$-packets}
As in~\S\ref{sub:bc_ps}, suppose $\Pi$ consists entirely of
constituents of the depth-zero principal series $\ind_B^G \lambda $.
Since each element of $\Pi$ has depth zero, $\ind_B^G \lambda $ and
hence $\lambda$ have depth zero by~\cite[Theorem~5.2]{moy-prasad2}.
It follows from~\cite{moy-prasad2} that for any $x\in\F$,
$(G_x,\lambda|_{M_0})$ is a $K$-type of each element of $\Pi$, where
$\finiteM$ is identified with $G_x/G_{x+}$.  Similarly,
$(\widetilde{G}_x, \tilde\lambda)|_{\widetilde{M}_0}$ is a $K$-type of
$\tilde{\pi} = \ind_{\widetilde{P}}^{\widetilde{G}}\tilde\lambda$ (see
Proposition~\ref{prop:bc_ps}), where
$\widetilde{G}_{x}/\widetilde{G}_{x+}$ is identified with
$\widetilde{\finiteM}$.  Denote by $\bar{\lambda}$ the character of
$\finiteM$ that inflates to $\lambda|_{M_0}$.  Then
$\tilde\lambda|_{\widetilde{M}_0}$ is the inflation to
$\widetilde{M}_0$ of the character $\tilde{\bar{\lambda}}$ of
$\widetilde{\finiteM}$.  As required, this is the Shintani lift of
$\bar{\lambda}$ from $\finiteM$ to $\widetilde{\finiteM}$.

\subsection{Singleton supercuspidal $L$-packets}
Now suppose that $\Pi$ is a singleton supercuspidal $L$-packet
$\{\pi\}$ of depth zero.  Then, by Proposition~\ref{prop:stable_sc},
$\pi$ is of the form $\ind_{G_y}^G\sigma$, where $\sigma$ is the
inflation to $G_y$ of a cubic cuspidal representation $\bar\sigma$ of
$\finiteG_y$.  Then $(G_y,\sigma)$ is a $K$-type of $\pi$
by~\cite[Prop.~6.2]{moy-prasad2}.  Similarly, it follows from
Proposition~\ref{prop:bc_stable} and~\cite[Prop.~6.2]{moy-prasad2}
that $(\widetilde{G}_y,\tilde{\sigma})$ is a $K$-type of
$\tilde{\pi}$, where $\tilde{\sigma}$ is the inflation to
$\widetilde{G}_y$ of the Shintani lift $\tilde{\bar\sigma}$ of
$\bar\sigma$ from $\finiteG_y$ to $\widetilde{\finiteG}_y$.  Hence the
theorem holds in this case.

\subsection{Supercuspidal $L$-packets of size four}
\label{sec:l-packet-4}
Recalling the notation of Proposition~\ref{prop:card4}, suppose that
$\Pi = \{\pi_0,\pi_1,\pi_2,\pi_3\}$ is a depth-zero supercuspidal
$L$-packet.  By~\cite[Prop.~6.2]{moy-prasad2},
$(G_y ,\sigma)$ and $(G_z,\sigma_i)$ are $K$-types for $\pi_0$ and the
$\pi_i$ ($i=1,2,3$), respectively.

According to Proposition~\ref{prop:bc_card4}, $\tilde{\pi}$ is the principal
series representation $\ind_{\widetilde{B}}^{\widetilde{G}}\chi^*$,
where $\chi^*$ is the inflation to $\widetilde{M}\cong
E^\times\times E^\times\times E^\times$ of
$$
\hat{\chi} = \tilde\chi_1\otimes
\tilde\chi_2 \otimes \tilde\chi_3
\in\Hom(\widetilde\finiteM,\C^\times).
$$
View $\chi^*|_{\widetilde{M}_0}$ as a character of
$\widetilde{G}_x$ (for any $x\in \F$) under the identification
$\widetilde{\finiteG}_x=\widetilde{\finiteM}$.
Then, by~\cite[Thm.~5.2]{moy-prasad2},
$(\widetilde{G}_x,\chi^*|_{\widetilde{M}_0})$ is a $K$-type for
$\tilde{\pi}$.
Since $\tilde{\pi}$ contains
$(\widetilde{G}_x,\chi^*|_{\widetilde{M}_0})$, it follows that, as a
representation of $\widetilde{\finiteB}_y$,
$\tilde{\pi}^{\widetilde{G}_{x+}}$ contains the character
$\hat{\chi}$ of $\widetilde{\finiteB}_y$.  Hence, by Frobenius
reciprocity, $\tilde{\pi}^{\widetilde{G}_{y+}}$ contains a
subrepresentation of
$\ind_{\widetilde{\finiteB}_y}^{\widetilde{\finiteG}_y}\hat{\chi}$.
But
$\ind_{\widetilde{\finiteB}_y}^{\widetilde{\finiteG}_y}\hat{\chi}$
is irreducible, as $\hat{\chi}$ is in general position,
so $\tilde{\pi}^{\widetilde{G}_{y+}}$ contains
$\ind_{\widetilde{\finiteB}_y}^{\widetilde{\finiteG}_y}\hat{\chi}$.
This is the Shintani lift of $\bar\sigma$ from $\finiteG_y$ to
$\widetilde{\finiteG}_y$ (see~\ref{sec:Shintani-finite}).  

Identify $\widetilde\finiteG_z$ with
$\widetilde\finiteH\subset\widetilde\finiteG$.  Now,
by~\ref{sec:Shintani-finite}, the Shintani lift of $\bar\sigma_i$ is
$\ind_{\widetilde{\finiteB}_z}^{\widetilde{\finiteG}_z}({}^w\hat{\chi})$
for an appropriate $w\in 
W_{k_F}(\underline{\widetilde{\finiteM}},\underline{\widetilde{\finiteG}})$.
The argument in the preceeding paragraph, applied to
$\ind_{\widetilde{B}}^{\widetilde{G}}(^w\chi^*)\cong
\ind_{\widetilde{B}}^{\widetilde{G}}\chi^*$ (where we identify
$W_{k_F}(\underline{\widetilde{\finiteM}},\underline{\widetilde{\finiteG}})$
and
$W_F(\underline{\widetilde{M}},\underline{\widetilde{G}})$), shows that
$\tilde{\pi}^{\widetilde{G}_{z+}}$ contains
$\ind_{\widetilde{\finiteB}_z}^{\widetilde{\finiteG}_z}({}^w\hat{\chi})$.

\subsection{$L$-packets and $A$-packets of size two}
Now suppose $\Pi $ is an $L$-packet of the form $\{ \pi^2(\lambda)
,\pi^s(\lambda)\}$ or an $A$-packet of the form $\{ \pi^n(\lambda)
,\pi^s(\lambda)\}$ for some $\lambda\in\Hom (M,\C^\times)$ of depth
zero (see case (\ref{sec:depth-zero-G}PS--2) and (\ref{sub:ps})).
Both $\pi^2(\lambda)$ and $\pi^n(\lambda)$ are constituents of the
principal series $\ind_B^G\lambda$.  It follows
from~\cite{moy-prasad2} that both $\ind_B^G\lambda$ and $\lambda$ have
depth zero and that for any $x$ in $\F$, $(G_x,\lambda|_{M_0})$ is a
$K$-type for both of these representations.  By
Proposition~\ref{prop:bc_card_2}, $\tilde{\pi}$ is always a
constituent of the principal series $
\ind_{\widetilde{B}}^{\widetilde{G}}\tilde\lambda.  $ Therefore, as
above, $(\widetilde{G}_{x},\tilde\lambda|_{\widetilde{M}_0})$ is a
$K$-type for $\tilde{\pi}$.  But $\tilde\lambda|_{\widetilde{M}_0}$ is
the inflation of $\tilde{\bar{\lambda}}\in\Hom
(\widetilde{\finiteM},\C^\times)$, where $\bar{\lambda}\in\Hom
(\finiteM,\C^\times)$ is the character that inflates to
$\lambda|_{M_0}$.  This shows that the theorem is true for
$\pi^2(\lambda)$ and $\pi^n(\lambda)$.

It remains to consider $\pi^s(\lambda)$ (both as an element of
$\{\pi^s(\lambda),\pi^2(\lambda)\}$ and as one of
$\{\pi^s(\lambda),\pi^n(\lambda)\}$).  Let $\lambda_1 ,\lambda_2$ be
the respective characters of $E^\times,E^1$ determined by $\lambda$
according to \eqref{eqn:lambda}.  Suppose first
that $\lambda_1|_{\OK_E^\times}$ is trivial.  Then
Proposition~\ref{prop:card2} implies that $(G_y,\sigma)$ is a $K$-type for
$\pi^s(\lambda)$, where $\sigma$ is the inflation to $G_y$ of
$\tau\cdot (\bar{\lambda}_2\circ\det)$.  Here
$\bar\lambda_2$ is the character of $k_E^1$ determined by 
$\lambda_2$, and $\tau$ is the
cuspidal unipotent representation of $\finiteG_y$.  From
\S\ref{sec:Shintani-finite}, the Shintani lift of $\tau\cdot
(\bar{\lambda}_2\circ\det_{\underline{\finiteG}_y})$ from $\finiteG_y$ to
$\widetilde{\finiteG}_y$ is
\begin{equation}
\label{eq:cusp_unip}
\tilde{\tau}\cdot
(\tilde{\bar{\lambda}}_2\circ\det_{\underline{\widetilde\finiteG}_y}),
\end{equation}
where 
$\tilde{\tau}$ is the unipotent representation of
$\widetilde{\finiteG}_y$ that is neither the trivial nor the Steinberg
representation.
Let $\tilde{\sigma}$ be the inflation of this representation to
$\widetilde{G}_y$.
Proposition~\ref{prop:bc_card_2} states that the base change lift
$\tilde{\pi}$ of
$\Pi$ is 
$$
\ind_{\widetilde{P}}^{\widetilde{G}}\tilde{\rho}',
$$
where
$\tilde{\rho}'$
is either a
one-dimensional representation $\tilde{\xi}'$ of $\widetilde{H}$ or
$\St_{\widetilde{H}}(\tilde{\xi}')$.

Suppose that $\tilde {\rho}' = \tilde{\xi}'$.  By
Proposition~\ref{prop:bc_card_2}, 
$$
\tilde{\xi}' = (\lambda_1 \tilde{\lambda}_2 |\cdot |_E^{\mp 1/2}
\circ\det_{\GL(2)} ) \otimes \tilde{\lambda}_2.
$$
Using Mackey's theorem and Frobenius reciprocity, we have
\begin{align}
\label{eq:hom}
\Hom_{\widetilde{G}_y}(\tilde{\sigma} ,\text{Res}_{\widetilde{G}_y}
\ind_{\widetilde{P}}^{\widetilde{G}}\tilde{\xi}')
& = 
\Hom_{\widetilde{G}_y}(\tilde{\sigma} ,
\ind_{\widetilde{P}\cap\widetilde{G}_y}^{\widetilde{G}_y}\tilde{\xi}')
\nonumber\\
& =  \Hom_{\widetilde{P}\cap\widetilde{G}_y}(\tilde{\sigma}
,\tilde{\xi}')\nonumber\\
& = 
\Hom_{\widetilde{P}\cap\widetilde{G}_y}
        (\tilde\sigma\cdot {\tilde\xi}^{\prime-1},{\bf 1}),
\end{align}
where we interpret $\tilde\sigma\cdot{\tilde\xi}^{\prime-1}$ as
the product of the restriction of each factor to
$\widetilde{P}\cap\widetilde{G}_y$.  Identify $\widetilde{\finiteG}_y$ 
with
$\widetilde{\finiteG}$.  Since $\lambda_1|_{\OK_E^\times}$ is trivial,
$\tilde{\xi}'|_{\widetilde{P}\cap\widetilde{G}_y}$ is the inflation to
$\widetilde{P}\cap\widetilde{G}_y$ of the character
$\tilde{\bar\lambda}_2\circ\det_{\widetilde{\underline\finiteH}}$ of
$\widetilde\finiteH$.
It follows that
$\tilde{\sigma}\cdot {\tilde\xi}^{\prime-1}$
is the restriction to $\widetilde{P}\cap\widetilde{G}_y$
of the
inflation to $\widetilde{G}_y$ of
$\tilde{\tau}$.
Since both
$\tilde{\sigma}\cdot {\tilde\xi}^{\prime-1}$ and ${\bf 1}$ are trivial on
$\widetilde{G}_{y+}$, \eqref{eq:hom} can be identified
with
$$
\Hom_{\widetilde{\finiteP}}(\tilde{\tau},{\bf 1}),
$$
where $\widetilde{\finiteP}$ is the parabolic subgroup of
$\widetilde{\finiteG}_y$ whose inverse image in $\widetilde{G}_y$
contains $\widetilde{P}\cap\widetilde{G}_y$.  By Frobenius
reciprocity,
$$
\Hom_{\widetilde{\finiteP}}(\tilde{\tau},{\bf 1})
= \Hom_{\widetilde{\finiteG}_y}
(\tilde{\tau},\ind_{\widetilde{\finiteP}}^{\widetilde{\finiteG}_y} {\bf 1}).
$$
It is easily seen that
$\ind_{\widetilde{\finiteP}}^{\widetilde{\finiteG}_y} {\bf 1}$ has two
irreducible components: the trivial representation and
$\tilde{\tau}$.  Hence 
$$
\dim_\C \Hom_{\widetilde{G}_y}(\tilde{\sigma} ,\text{Res}_{\widetilde{G}_y}
\ind_{\widetilde{P}}^{\widetilde{G}}\tilde{\xi}') = 1.
$$
In particular, as a representation of
$\widetilde{G}_y$, $\tilde{\pi}^{\widetilde{G}_{y+}}$
must contain $\tilde{\sigma}$, as required.

Now suppose that $\tilde{\rho}' =
\St_{\widetilde{H}}(\tilde{\xi}')$.
By
Proposition~\ref{prop:bc_card_2}, the representations
$\ind_{\widetilde{P}}^{\widetilde{G}}\tilde{\rho}'$ and
$\ind_{\widetilde{P}}^{\widetilde{G}}\tilde{\xi}'$ are the irreducible
constituents of
$\ind_{\widetilde{B}}^{\widetilde{G}}\tilde\lambda$.
For all $x\in\F$,
$$
\text{Res}_{\widetilde{G}_y}
        \ind_{\widetilde{B}}^{\widetilde{G}} \tilde\lambda
=
\ind_{\widetilde{B}\cap \widetilde{G}_y} ^{\widetilde{G}_y} \tilde\lambda
=
\ind_{\widetilde{G}_x}^{\widetilde{G}_y}
\ind_{\widetilde{B}\cap \widetilde{G}_y} ^{\widetilde{G}_x}
\tilde\lambda,
$$
which contains $ \ind_{\widetilde{G}_x} ^{\widetilde{G}_y}
\tilde\lambda $, the inflation to $\widetilde{G}_y$ of the
representation $\ind_{\widetilde{\finiteB}_y}^{\widetilde{\finiteG}_y}
\tilde{\bar\lambda}$.  (Here $\bar\lambda$ is the character of
$\finiteM$ determined by $\lambda$.)  Moreover,
$\ind_{\widetilde{\finiteB}_y}^{\widetilde{\finiteG}_y}
\tilde{\bar{\lambda}}$ contains two copies of the representation
(\ref{eq:cusp_unip}) since 
$$
\bigl(\ind_{\widetilde{\finiteB}_y}^{\widetilde{\finiteG}_y}
  \tilde{\bar{\lambda}}\bigr)\cdot \bigl(\tilde{\bar\lambda}_2^{-1}
  \circ\det_{\widetilde{\underline\finiteG}_y}\bigr) =
\ind_{\widetilde{\finiteB}_y}^{\widetilde{\finiteG}_y}{\bf 1}
$$
contains two copies of $\tilde{\tau}$.  Therefore,
$$
\dim_\C\Hom_{\widetilde{G}_y}
(\tilde\sigma,\ind_{\widetilde{B}}^{\widetilde{G}}\tilde\lambda) \geq 2.
$$
Since
$$
\dim_\C\Hom_{\widetilde{G}_y}(\tilde{\sigma} ,\text{Res}_{\widetilde{G}_y}
\ind_{\widetilde{P}}^{\widetilde{G}}\tilde{\xi}') = 1,
$$
it follows that 
$$
\Hom_{\widetilde{G}_y}(\tilde{\sigma} ,\text{Res}_{\widetilde{G}_y}
\ind_{\widetilde{P}}^{\widetilde{G}}\tilde{\rho}')\neq 0.
$$
Hence, as above, $\tilde{\pi}^{\widetilde{G}_{y+}}$
must contain $\tilde{\sigma}$.

On the other hand, suppose that $\lambda_1|_{\OK_E^\times}$ is not
trivial.  Let $\bar\lambda_1'$ be the character of $k_E^1$ determined
by $\bar{\lambda}'_1\circ\N = \bar{\lambda}_1$.  Let $\sigma$ be the
inflation to $G_z$ of the cuspidal representation $\bar\sigma$ of
$\finiteG_z$ with character $-R_\finiteC^{\finiteG_z}\,\chi$, where
$$
\chi = \bar{\lambda}'_1\bar{\lambda}_2 \otimes
\bar{\lambda}'_1\bar{\lambda}_2 \otimes \bar{\lambda}_2.
$$
Then $(G_z,\sigma)$ is a
$K$-type for $\pi^s(\lambda)$.
By \S\ref{sec:Shintani-finite}, the Shintani lift $\tilde{\bar\sigma}$
of $\bar\sigma$ from $\finiteG_z$ to $\widetilde{\finiteG}_z$ is
$
\ind_{\widetilde{\finiteB}_z}^{\widetilde{\finiteG}_z}
  \hat{\chi}$,
where
$
\hat{\chi} = \bar{\lambda}_1\tilde{\bar{\lambda}}_2 \otimes
\bar{\lambda}_1\tilde{\bar{\lambda}}_2 \otimes
\tilde{\bar{\lambda}}_2$.

Now in both the $L$-packet and $A$-packet cases, $\tilde{\pi}$ is a
constituent of $\ind_{\widetilde{B}}^{\widetilde{G}}\tilde\lambda$ by
the proof of~\cite[Lemma 12.7.6]{rogawski}.  Hence
$(\widetilde{G}_x,{}^w\tilde\lambda|_{\widetilde{M}_0})$ is a $K$-type
of $\tilde{\pi}$, where $x$ is any point in $\F$ and $w\in
W_F(\underline{\widetilde{M}},\underline{\widetilde{G}})$.
As in the case above where $\pi$
is supercuspidal of type \depthzerosc{4}, it follows by Frobenius
reciprocity that $\tilde{\pi}^{\widetilde{G}_{z+}}$ contains a
subrepresentation of $
\ind_{\widetilde{\finiteB}_z}^{\widetilde{\finiteG}_z}
({}^w\tilde{\bar{\lambda}}), $ since
$^w\tilde\lambda|_{\widetilde{M}_0}$ is the inflation of the character
$^w\tilde{\bar{\lambda}}$ of $\widetilde\finiteM$ (where we have
identified $W_F(\underline{\widetilde{M}},\underline{\widetilde{G}})$ and
$W_{k_F}(\underline{\widetilde{\finiteM}},\underline{\widetilde{\finiteG}})$).
Using the fact
that $\lambda_1|_{\OK_F^{\times}}$ is trivial, one finds that
$$
\tilde{\bar{\lambda}} =
\bar{\lambda}_1 \tilde{\bar{\lambda}}_2
\otimes
\tilde{\bar{\lambda}}_2\otimes
\bar{\lambda}_1 \tilde{\bar{\lambda}}_2
$$
so for an appropriate $w$,
$$
{}^w\tilde{\bar{\lambda}} =
\bar{\lambda}_1 \tilde{\bar{\lambda}}_2
\otimes
\bar{\lambda}_1 \tilde{\bar{\lambda}}_2\otimes
\tilde{\bar{\lambda}}_2.
$$
Thus
$
\ind_{\widetilde{\finiteB}_z}^{\widetilde{\finiteG}_z}
({}^w\tilde{\bar{\lambda}}) = \tilde{\bar\sigma},
$
so 
$\tilde{\pi}^{\widetilde{G}_{y+}}$ contains the irreducible
representation $\tilde{\bar\sigma}$.

\section{On Induced characters of nonconnected groups}
\label{sec:induced}

Let $G$ now denote the group of rational points of a reductive group
defined over a nonarchimedean local field.
In particular, we do not assume that $G$ is connected.
We will, however, assume that $G$ is a semidirect product
of its connected component $G^0$ and its component group $\Gamma$,
and that $G$ has a $\Gamma$-invariant special parahoric subgroup.
Suppose $H$ is an open subgroup of $G$ that is compact modulo the center
of $G$.  Let $\rho$ denote an irreducible, smooth representation
of $H$, and let $\pi$ denote the compactly induced representation
$\ind_H^G\rho$ of $G$.  
Let $K$ denote a compact open subgroup of $G$.
\begin{prop}
\label{induced_char}
For $g\in G\reg$,
$$
\theta_\pi(g) =
\sum_{a\in K \backslash  G / H}
\bigl(\sum_{b \in KaH/H} \theta_\rho(b^{-1}gb)\bigr) ,
$$
where $\theta_\rho$ is extended to $G$ by zero.  For each $g$, all
but finitely many terms of the inner sum vanish.
\end{prop}

\begin{proof}
This is identical to the proof of Theorem~A.14
of~\cite{bushnell-henniart:local-tame-1}.
\end{proof}


\begin{thebibliography}{10}
\expandafter\ifx\csname url\endcsname\relax
  \def\url#1{\texttt{#1}}\fi
\expandafter\ifx\csname urlprefix\endcsname\relax\def\urlprefix{URL }\fi

\bibitem{adler-debacker:mk-theory}
J.~D. Adler, S.~DeBacker, Murnaghan-{K}irillov theory for supercuspidal
  representations of tame general linear groups, J. Reine Angew. Math. 575
  (2004) 1--35.

\bibitem{badulescu:thesis}
A.~Badulescu, Correspondance entre $\hbox{GL}_n$ et ses formes int\'erieures en
  caract\'eristique positive, Ph.D. thesis, Universit\'e de Paris-Sud (1999).

\bibitem{bruhat-tits:one}
F.~Bruhat, J.~Tits, Groupes r\'eductifs sur un corps local {I}: {D}onn\'ees
  radicielles valu\'ees, Publ.~Math.~IHES 41 (1972) 5--251.

\bibitem{bruhat-tits:two}
F.~Bruhat, J.~Tits, Groupes r\'eductifs sur un corps local {II}: Sch\'emas en
  groupes. {E}xistence d'une donn\'ee radicielle valu\'ee, Publ.~Math.~IHES 60
  (1984) 197--376.

\bibitem{bushnell-henniart:local-tame-1}
C.~J. Bushnell, G.~Henniart, Local tame lifting for $\hbox{GL}({N})$. {I}.
  {S}imple characters, Inst. Hautes \'Etudes Sci. Publ. Math. (1996) 105--233.

\bibitem{bushnell-henniart:local-tame-2}
C.~J. Bushnell, G.~Henniart, Local tame lifting for $\hbox{GL}(n)$. {II}.
  {W}ildly ramified supercuspidals, Ast\'erisque 254 (1999).

\bibitem{bushnell-henniart:bc-glp}
C.~J. Bushnell, G.~Henniart, Explicit unramified base change: {${\rm GL}(p)$}
  of a {$p$}-adic field, J. Number Theory 99~(1) (2003) 74--89.

\bibitem{bushnell-henniart:local-tame-4}
C.~J. Bushnell, G.~Henniart, Local tame lifting for {${\rm GL}(n)$}. {IV}.
  {S}imple characters and base change, Proc. London Math. Soc. (3) 87~(2)
  (2003) 337--362.

\bibitem{bushnell-kutzko:smooth}
C.~J. Bushnell, P.~C. Kutzko, Smooth representations of reductive $p$-adic
  groups: structure theory via types, Proc. London Math. Soc. (3) 77~(3) (1998)
  582--634.

\bibitem{debacker}
S.~DeBacker, Some applications of {B}ruhat-{T}its theory to harmonic analysis
  on a reductive {$p$}-adic group, Michigan Math. J. 50~(2) (2002) 241--261.

\bibitem{dkv:amt}
P.~Deligne, D.~Kazhdan, M.-F. Vign{\'e}ras, Repr\'esentations des alg\`ebres
  centrales simples $p$-adiques, in: Representations of reductive groups over a
  local field, Hermann, Paris, 1984, pp. 33--117.

\bibitem{deligne-lusztig}
P.~Deligne, G.~Lusztig, Representations of reductive groups over finite fields,
  Ann.~of Math. 103 (1976) 103--161.

\bibitem{digne0}
F.~Digne, Shintani descent and {$\mathcal{L}$} functions on {D}eligne-lusztig
  varieties, in: P.~Fong (Ed.), The {A}rcata Conference on Representations of
  Finite Groups, Vol. 47, {P}art 2 of Proc.~Symp.~Pure Math., Amer.~Math.~Soc.,
  Providence, RI, 1987, pp. 61--68.

\bibitem{digne}
F.~Digne, Descente de {S}hintani et restriction des scalaires, J.~London
  Math.~Soc. 59~(2) (1999) 867--880.

\bibitem{ennola}
V.~Ennola, On the characters of the finite unitary groups, Ann.\ Acad.\ Sci.\
  Fenn.\ Ser.\ A I 323 (1963).

\bibitem{gyoja}
A.~Gyoja, Liftings of irreducible characters of finite reductive groups, Osaka
  J.~Math 16 (1979) 1--30.

\bibitem{kawanaka}
N.~Kawanaka, On the irreducible characters of the finite unitary groups, J.
  Math. Soc. Japan 29~(3) (1977) 425--450.

\bibitem{Keys}
C.~D. Keys, Principal series representations of special unitary groups over
  local fields, Compositio Math. 51~(1) (1984) 115--130.

\bibitem{kim-ps:theta-10}
J.-L. Kim, I.~I. Piatetski-Shapiro, Quadratic base change of {$\theta\sb
  {10}$}, Israel J. Math. 123 (2001) 317--340.

\bibitem{morris:intertwining}
L.~Morris, Tamely ramified intertwining algebras, Invent. Math. 114 (1993)
  1--54.

\bibitem{moy-prasad2}
A.~Moy, G.~Prasad, Jacquet functors and unrefined minimal {$K$}-types,
  Comment.\ Math.\ Helv. 71~(1) (1996) 98--121.

\bibitem{rogawski:amt}
J.~Rogawski, Representations of $\hbox{GL}(n)$ and division algebras over a
  $p$-adic field, Duke Math. J. 50~(1) (1983) 161--196.

\bibitem{rogawski}
J.~Rogawski, Automorphic Representations of Unitary Groups in Three Variables,
  Vol. 123 of Annals of Math.\ Studies, Princeton University Press, Princeton,
  New Jersey, 1990.

\bibitem{schneider-stuhler}
P.~Schneider, U.~Stuhler, Representation theory and sheaves on the
  {B}ruhat-{T}its building, Publ. Math. IHES 85 (1997) 97--191.

\bibitem{shintani}
T.~Shintani, Two remarks on irreducible characters of finite general linear
  groups, J.\ Math.\ Soc.\ Japan 28~(2) (1976) 396--414.

\bibitem{silberger-zink:level-zero-matching}
A.~J. Silberger, E.-W. Zink, An explicit matching theorem for level zero
  discrete series of unit groups of $p$-adic simple algebras, preprint, 2004.

\bibitem{srinivasan}
B.~Srinivasan, Representations of Finite {C}hevalley Groups, Vol. 764 of
  Lecture Notes in Mathematics, Springer-Verlag, Berlin, 1979.

\bibitem{waldspurger:fourier}
J.-L. Waldspurger, Transformation de {F}ourier et endoscopie, J. Lie Theory
  10~(1) (2000) 195--206.

\bibitem{waldspurger:nilpotent}
J.-L. Waldspurger, Int\'egrales orbitales nilpotentes et endoscopie pour les
  groupes classiques non ramifi\'es, Ast\'erisque 269 (2001).

\end{thebibliography}
\end{document}